\documentclass[a4paper,twoside,11pt]{amsart}
\usepackage{multicol}
\usepackage[all]{xy}
\usepackage{appendix}
\usepackage{graphicx}
\usepackage{tikz} 
\usepackage{tikz-3dplot}

\usepackage[pagebackref]{hyperref}
\hypersetup{
	colorlinks=true,
	linkcolor=cyan,
	citecolor=cyan
			}

\usepackage[section]{placeins} % Prevents placing floats before or after the section.

\usepackage[makeindex]{imakeidx}

\usepackage{amsfonts, amsmath, amssymb, amsthm, amsopn, latexsym, graphicx, mathrsfs,  verbatim, enumerate, url, stmaryrd, enumitem, tikz,mathtools} % Standard packages
\usepackage[normalem]{ulem}
\usepackage[all]{hypcap} %So that hyperlinks to figures go to the figure and not the caption
\usepackage{pst-all} %Package to draw figures

\usepackage[pagewise]{lineno} %To number all lines
\usepackage{color}

\usepackage{booktabs} %Tablas profesionales

\usepackage{caption}
\usepackage{pdfpages} %Include pdf or pages in a tex document 
\definecolor{green}{rgb}{0.0, 0.5, 0.0}
\definecolor{yellow}{rgb}{0.8, 0.33, 0.0}
\definecolor{purple}{RGB}{172, 51, 255}
 \usepackage[a4paper]{geometry}

\makeatletter
\def\subsection{\@startsection{subsection}{2}%
 \z@{.5\linespacing\@plus.7\linespacing}{.3\linespacing}%
 {\normalfont\bfseries}}
\makeatother

 \usepackage[active]{srcltx}

\usepackage{amsthm,amsmath,amssymb,amscd}
\usepackage{graphicx}
\usepackage{epic}

\newtheorem{theorem}{Theorem}[section]
\newtheorem{proposition}[theorem]{Proposition}
\newtheorem{lemma}[theorem]{Lemma}
\newtheorem{corollary}[theorem]{Corollary}

\theoremstyle{definition}
\newtheorem{definition}[theorem]{Definition}
\newtheorem{notation}[theorem]{Notation}

\newtheorem{example}[theorem]{Example}
\theoremstyle{remark}
\newtheorem{remark}[theorem]{Remark}

\newtheorem{thmno}{Theorem}

\numberwithin{equation}{section}

\def\1{\mathbf{1}}

\def\a{\alpha}

\def\b{\beta}

\def\C{{\mathbf C}}

\def\D{\Delta}
\def\d{\delta}

\def\g{\gamma}

\def\l{\lambda}

\def\N{{\mathbf N}}
\def\Newton{{\mathcal N}}

\def\orb{\mathrm{orb}}

\def\orden{\mathrm {ord}}

\def\p{\pi}

\def\Q{{\mathbf Q}}

\def\r{\rho}

\def\R{{\mathbf R}}

\def\x{x}

\def\z{\zeta}
\def\\X{D}

\def\Z{{\mathbf Z}}

\def\s{\sigma}
\def\t{\tau}

\def\elem(#1,#2){  \{ \frac{#1}{\overline{\ #2\ }}\}  }

\def\C {\mathbb C}
\def\N {\mathbb N}
\def\R {\mathbb R}

\def\Q {\mathbb Q}

\def\Z {\mathbb Z}

\def\a{\alpha}

\def\b{\beta}

\def\D{\Delta}
\def\d{\delta}

\def\g{\gamma}

\def\l{\lambda}

\def\Newton{\mathcal N}

\def\X{D}

\def\p{\pi}

\def\r{\rho}

\def\s{\sigma}
\def\t{\tau}
\newcommand{\trop}{\mathrm{trop}}

\def\z{\zeta}

\makeindex
\def\elem(#1,#2){  \left\{ \frac{#1}  {\overline {\ #2\ }} \right\} }

\title{Quasi-ordinary hypersurfaces, multiplier ideals and local tropicalizations}
\author{Pedro D. Gonz\'alez P\'erez}
\address{Instituto de Matematica Interdisciplinar y Departamento de \'Algebra, Geometr\'{\i}a y Topolog\'{\i}a, Facultad de Ciencias Matem\'aticas, Universidad Complutense de Madrid, Plaza de las Ciencias 3, Madrid 28040, Espa\~na.}

\email{pgonzalez@mat.ucm.es}

\author[M. Robredo Buces] {Miguel Robredo Buces}
\address{Instituto de Ciencias Matem\'aticas (CSIC-UAM-UC3M-UCM), Calle Nicol\'as Cabrera,  13-15
Campus de Cantoblanco, 
28049 Madrid Spain}
\email{miguel.robredo@icmat.es}

\thanks{The authors are grateful to the referee for the detailed revision of a previous version of the paper.
 The first author is supported by the Spanish grant of MCIN with reference 
(PID2020-114750GB-C32/AEI/10.13039/501100011033). The second author was supported by the grant  SEV-2015-0554}

\keywords{multiplier ideals, jumping numbers, quasi-ordinary hypersurface singularities, toroidal resolutions, local tropicalization}

\date{\today}

\subjclass[2000]{32S25, 32S45, 32S35, 14M25}

\begin{document}
% \linenumbers

\begin{abstract}
In this paper we describe the multiplier ideals and  jumping numbers 
associated with an irreducible germ of quasi-ordinary hypersurface $(\X, 0) \subset (\C^{d+1}, 0)$
by using a toroidal embedded resolution. The approach is motivated by  Howald's description of the multiplier ideals of monomial ideals. 
We show that the multiplier ideals of $D$ can be expressed in terms of a finite sequence of Newton polyhedra associated with the total transform of $\X$ in the toroidal resolution process.
  We prove that the multiplier ideals are  \textit{generalized monomial ideals} with respect to a complete sequence of semi-roots. This is a finite sequence of functions which determines a system of generators of the semigroup of the quasi-ordinary hypersurface. 
We express these results in terms of the local tropicalization associated with the embedding 
of $\C^{d+1}$ 
defined by this sequence. We prove that the local tropicalization is the support of a 
fan of the lattice $\Z^{d+g+1}$, which is determined by the embedded topological type of 
$(\X, 0) \subset (\C^{d+1}, 0)$.
\end{abstract}
\maketitle

\bigskip
{\em \hfill This paper is dedicated to Bernard Teissier. \smallskip}

\tableofcontents

\section* {Introduction}

Let $Y$ be a smooth complex algebraic variety and $D$ be an integral reduced effective divisor on $Y$.
 Let 
 $ \Psi :Z\rightarrow Y$
be a  proper birational map,  with $Z$ a smooth variety. 
The map $\Psi$ is an \textit{embedded resolution} of $D$ if 
\begin{enumerate}[label=(\alph*)]
\item $\Psi^{*} \mathcal{O} (-D) =\mathcal{O}_{Z}(-F)$, where 
$F$ is a simple normal crossing divisor,
\item \label{c-a} $\Psi$ defines an isomorphim $Z \setminus F \rightarrow Y \setminus D$,
\item \label{c-b}
the restriction of $\Psi$ to the strict transform  of $D$ is an isomorphism over the smooth locus of $D$.
\end{enumerate}
 Notice that if 
 $\Psi$ is 
an embedded resolution of $D$,  then the prime  components of the exceptional divisor  of $\Psi$ are 
components of $F$. There is a also stronger version of embedded resolution where one asks that
 the exceptional divisor of $\Psi$ maps to the singular locus of $D$, instead of 
conditions \ref{c-a} and \ref{c-b}.
The existence of embedded resolutions follows from Hironaka's work on resolution of singularities.
We refer to the Spivakovsky's paper \cite{S20} for an introduction to  resolution of singularities.

\medskip 

The  \textit{multiplier ideal} $\mathcal{J}(\xi D)_o$ associated with $D$ and $\xi \in \Q_{\geq 0}$ at the point $o \in Y$ can be expressed by  valuative conditions in terms of a fixed embedded resolution $\Psi$. It consists of functions  
 $h $ in the ring $  \mathcal{O}_{Y, o}$ of germs of holomorphic functions of $Y$ at $o$ such that 
\begin{equation} \label{eq: val-ine}
  \nu_{E_i} (h)  + \lambda_{E_i}   > \xi \nu_{E_i} (f)
\end{equation}
where $f$ is a defining function of $D$ at $o$, $E_i$ runs through the prime components 
of $F$, $ \nu_{E_i}$ denotes the associated vanishing order valuation, 
$\lambda_{E_i} := 1+ \nu_{E_i}  (K_{\Psi})$ is the \textit{log-discrepancy} of $E_i$, 
and   $K_{\Psi}$ is the relative canonical divisor. 
The study of multiplier ideals has become a central aspect of birational geometry, since 
they provide a subtle measure of the singularities of the pair $( Y, D)$.
The multiplier ideals can be 
associated more generally to an ideal sheaf on $\mathcal{O}_{Y}$ even when $Y$ has normal singularities
(see \cite{dFH09,  BdFF15}).
 See \cite{BL04} for an introduction to the subject, Lazarsfeld's book \cite{PAG} for general reference and 
 \cite{B12} for the relations with other topics in singularity theory.

\medskip 

There  exists an increasing sequence $(\xi_{i})_i$ of positive rational numbers, called 
the \textit{jumping numbers} of $D$ at $o$, 
such that if  $\xi_i \leq \xi <   \xi_{i+1}$ then
$
\mathcal{J} (\xi D )_o = 
	\mathcal{J} (\xi_{i} D )_o
	\supsetneq \mathcal{J}( \xi_{i+1} D)_o$ (see \cite{ELSVJumpingCoefficients}).
 If $Y$ is a complex surface there are many results about multiplier ideals and
 their associated jumping numbers (see for instance  \cite{JarJNSCI, TuckerJNRatSing, GM10, AADG17, JNN, AAD, Blanco,FJ05, GGGR24}).
 
 \medskip

This paper describes the multiplier ideals and jumping numbers  associated with an irreducible quasi-ordinary hypersurface singularity. It is based on the PhD thesis of the second author \cite{R19}.

\medskip

An equidimensional germ  $(\X,o)$ of a complex analytic variety  of
dimension $d$ is \textit {quasi-ordinary} (q.o.) if there exists a
finite map $\pi: (\X,o)\rightarrow(\C^d,0)$ which is unramified
outside a normal crossing divisor of $(\C^d,0)$.
Q.o.~surface singularities originated in Jung's  classical method to parametrize and resolve
surface singularities (see \cite{Jung}).  
In  what follows we assume  that $(\X,o) \subset (\C^{d+1},0)$ is an irreducible quasi-ordinary hypersurface. By the Abhyankar-Jung Theorem,  there exists local coordinates $(x, y):= (x_1,\dots, x_d, y)$ at the origin of $\C^{d+1}$ such 
that $\X$ admits a fractional power series parametrization of the form
$y = \zeta(x^{1/n})$ for some integer $n$ where $x^{1/n} =  ({x_1}^{1/n}, \dots, x_{d}^{1/n})$, 
see \cite{A55, PR12}. 
The series $\zeta:= \zeta(x^{1/n})$, which is called a \textit{quasi-ordinary branch}, has a $g \geq 0$ 
 \textit{characteristic monomials}, which generalize the notion of Puiseux  fractional characteristic exponents
 in the case of plane branches (see \cite{LipmanTopInvQO}).  
The characteristic  monomials of a \textit{normalized} q.o. branch 
 classify the embedded topology of the germ $(  D , 0) \subset
 (\C^{d+1}, 0)$ (see \cite{LipmanTopInvQO,Gau}). 
 We may assume that $\z(x)$ is \textit{normalized} by changing the projection $\pi$ (see \cite{Gau}).

\medskip 

We study the quasi-ordinary hypersurface $D$ by using a \textit{complete sequence of semi-roots}. This is a sequence of functions $x_1, \dots,  x_{d+g+1} \in \C\{ x, y \}$, defining hypersurfaces $D_1, \dots D_{d+g+1}$ such that $D_1, \dots, D_d$ the coordinate hyperplanes and $D_{d+1}, \dots, D_{d+g+1} = D$ are certain quasi-ordinary hypersurfaces.
This sequence generalizes the maximal contact curves of higher genus in the case of plane branches (see \cite{Lejeune, PP03}). 
This sequence determines generators of the \textit{semigroup} associated to the q.o. hypersurface.
It generalizes  in a subtler way, the classical notion of semigroup of a plane branch (see \cite{SgrPG, PPPDuke, KiMiSgr}). 
In addition,  $(x_1, \dots, x_{d+1})$ is a system local coordinates 
at the origin of $(\C^{d+1},0)$, and  this sequence provides suitable local coordinates 
at certain points appearing in the  toroidal embedded resolution process of the q.o. hypersurface.

\medskip 

Let us outline  the construction of  the toroidal embedded resolutions of $\X$ (see Section \ref{resolution}, and \cite{GPRQo, MMFGG}).
 We define a sequence of pairs  $(Z_\ell, o_\ell)$, where $Z_\ell$ is a normal variety of dimension $d+1$ and $o_\ell \in Z_\ell$ is a closed point, 
 together with modifications 
 \[
\psi_\ell: Z_\ell \longrightarrow Z_{\ell-1}, 
 \]
 for $\ell = 0, \dots, g$, where $Z_{0} = \C^{d+1}$. 
 The modification $\psi_1$ is the \textit{Newton modification} associated 
 with a defining function of $\X$ with respect to the local system of coordinates $(x_1, \dots, x_{d+1})$, that is, 
 $\psi_1$ is  the modification determined by the dual fan of the Newton polyhedron of  $x_{d+g+1}$ with respect to these coordinates.
By using the sequence $ \x_1, \dots, \x_{d+g+1} $ we show that the local ring 
 of germs of analytic functions at $o_\ell \in Z_\ell$ is isomorphic 
 the local ring   $R_\ell$  of a normal toric variety of dimension $d+1$ which is determined by the characteristic monomials,  for $\ell \in \{1, \dots, g \}$. 
 This provides a local toric structure of $Z_\ell$ at the point $o_\ell$.
If  $h \in \C \{ x_1, \dots, x_d, x_{d+1}  \}$, we denote by  $\mathcal{N}_{\ell} (h)$ the Newton polyhedron
   of a power series  defining the total transform of $h$ in the toric ring $R_\ell$.
The modification $\psi_\ell$ is the Newton modification associated with a defining function of 
the total transform of $\X$ on $Z_{\ell-1}$ with respect to the toric structure defined by  $R_{\ell -1}$.
The composition
\[
  \Psi_{g} := \psi_{1} \circ \cdots \circ \psi_{g} : \, Z_g \to Z_0
\]
is a toroidal embedded normalization of the quasi-ordinary hypersurface singularity $\X$ (see \cite{GPRQo}).
\medskip

Let us denote by $D_i$ the hypersurface defined by $x_i = 0$, for $i \in \{ 1, \dots, d+g +1\}$.
We have a  \textit{toroidal embedding without self-intersection} associated to  the toroidal embedded normalization,  whose \textit {boundary divisor} $\partial Z_g$ is  the reduced divisor of the total transform of $D_1 +\cdots + D_{d+g+1}$.
The variety $Z_g$ has only  toric singularities. A toroidal modification $\Psi_{\mathrm{reg}}$ 
defined by a regular subdivision of the conic polyhedral complex $\bar{\Theta}$ of this toroidal embedding is a resolution of the singularities of $Z_g$. The modification $\Psi = \Psi_{\mathrm{reg}} \circ \Psi_g$ is an embedded resolution of the quasi-ordinary hypersurface $\X$ (see \cite{GPRQo}).

\medskip

We define a set of \textit{quasi-monomial valuations} associated with the toroidal embedded normalization. We prove that $x_1, \dots, x_{d+g+1}$ is a generating sequence for these quasi-monomial valuations (see Section \ref{sec-qmv}).
In Section \ref{Sect:log-discrepancy vector}, we describe the log-discrepancies of the components of the boundary divisor $\partial Z$  in terms of the characteristic monomials (see Section \ref{Sect:log-discrepancy vector}). More generally, we find a monomial $x^{\lambda_\ell} \in R_\ell$ which determines 
the log-discrepancies of the quasi-monomial divisorial valuations which appear at the step $\ell$-th of the toroidal resolution (see Proposition \ref{prop:log-discrepancy vector}). 
This is a generalization of a result of Favre and Jonsson  in the plane curve case \cite[Prop. D.1]{TVT}, see also \cite{GGP25, GGP19}.

\medskip 
The main contributions about the multiplier ideals are given in Section \ref{sec: miqh}. 
The following theorem (see Theorem \ref{Thm: reduction}) is inspired  by Howald's results \cite{MMI, Howald:Non-degenerate} about multiplier ideals of monomial ideals and 
 of Newton non-degenerate hypersurface singularities. 
 \begin{thmno} 
Take $\xi \in \Q_{\geq 0}$ and $h \in  \mathcal{O}_{\C^{d+1},0}$.
The following conditions are equivalent: 
\begin{enumerate}[label=(\alph*)]
\item  \label{a-0}
$h \in \mathcal{J}(\xi \X)_0$.
\item \label{b-0}
The polyhedron $ \Newton_{\ell -1} (h) + \lambda_{\ell -1}$ is contained in the interior of  the polyhedron  $ \xi \Newton_{\ell -1}( x_{d+g+1} )$, for 
$\ell \in \{1, \dots, g+1\}$.
\item   \label{c-0}
The valuative inequality \eqref{eq: val-ine}
holds for any prime divisor
 $E_i$ in the support of the boundary divisor $\partial Z_g$ of the toroidal embedded normalization. \end{enumerate}
\end{thmno} 
This theorem shows 
 that many of exceptional divisors of the toroidal embedded resolution $\Psi$ are not relevant in the definition of  the multiplier ideals.  
This is a higher dimensional generalization of Theorem 4.9 and Corollary 4.17 of our joint 
paper with Gonz\'alez Villa and Guzm\'an Dur\'an
\cite{GGGR24}. This result is inspired by 
a theorem of  Smith and Thompson, which states that an exceptional prime divisor 
 $E$ of 
 the minimal embedded resolution of a plane curve singularity $C$ \textit{contributes} to the sequence of jumping numbers if and only if $E$ intersects at least three other components of the total transform of $C$ (see \cite{STIrrelExcDiv}). One can compare this result with the description of a finite set of \textit{candidate poles} of the topological zeta function of   a defining function of the  quasi-ordinary hypersurface $D$ at $0$, which are determined by certain components of the boundary divisor $\partial Z_g$, see \cite[Cor. 4.18]{MMFGG}.

\medskip 

A \textit{generalized monomial in the semi-roots}  is of the form $\mathcal{M} = x_1^{i_1} \dots x_{d+g+1}^{i_{d+g+1}}$,  for $i_1, \dots, i_{d+g+1} \in \Z_{\geq 0}$. 
Note that $\mathcal{M}$  \textit{is not a monomial} of  the ring of power series $\C\{ x_1, \dots, x_{d+1} \} $ in general.
We associate to $\mathcal{M}$ the number:
\[
\xi_{\mathcal{M}}: = \min \left\{ ({\nu_{E_i}} ( \mathcal{M} ) + \lambda_{E_i})({  \nu_{E_i} ( x_{d+g+1} ) })^{-1}  \mid E_i \mbox{ is an exceptional component of } \partial Z_g 
 \right\}.
 \]
We prove the following description of the generators  of the multiplier ideals of $\X$ at $o$ and of their jumping numbers (see Theorem \ref{rem: mon-ideal}). It is obtained by 
 generalizing the approach of   \cite[Th. 4.20]{GGGR24}  in the case of plane curve singularities.

\begin{thmno} Take a rational number $\xi \in (0,1)$.  
\begin{enumerate}[label=(\alph*)]
 \item   
 The multiplier ideal  $\mathcal{J}(\xi \X)_0$ is generated by those generalized monomials 
 $\mathcal{M}$ in $x_1, \dots, x_{d+g+1}$ such that  
 $\xi < \xi_{\mathcal{M}} \leq  \xi + 1$. 
 \item  
 A rational number $\xi$ is a jumping number of the multiplier ideals of $\X$ at $0$ if and only if there exist 
 a generalized monomial $\mathcal{M}_0$ in $x_1, \dots, x_{d+g+1}$ such that  $\xi = \xi_{\mathcal{M}_0}$. 
\end{enumerate}
\end{thmno} 

As a consequence of this result the computation 
of the multiplier ideals and the associated jumping numbers of $\X$ 
boils down to an optimization problem in terms of log discrepancies 
of the components of the boundary divisor $\partial Z_g$ of the toroidal embedded normalization 
and the values of the corresponding exceptional divisors on the functions $x_1, \dots, x_{d+g+1}$.

\medskip
In Section \ref{sec: loctrop}, we apply to this case the theory of local tropicalization developed by Popescu-Pampu and Stepanov in  \cite{PS13, PS25}.  
The local tropicalization was used in \cite{dFGM23} to show the existence of
embedded resolutions of curve singularities with one toric morphism, after a suitable reembedding, generalizing a theorem of Goldin and Teissier in the case of plane branches \cite{GT}. 
Cueto, Popescu-Pampu and Stepanov gave another proof of the main result of \cite{dFGM23},  as a corollary of their description of the local tropicalizations of  some surface singularities \cite{CPS24}. 
\medskip 

We consider the embedding of $(\C^{d+1},0) \hookrightarrow (Y, 0) \subset (\C^{d+g+1},0)$ 
defined by the complete sequence of semi-roots $(x_1, \dots, x_{d+g+1})$. This embedding is inspired by 
the works of Teissier on plane branches, valuations,  and local uniformization (see
\cite{GT, T14}).
The first author proved that there there is one toric modification of $\C^{d+g+1}$, which induces an 
embedded resolution of the quasi-ordinary hypersurface $(\X,o) \hookrightarrow (Y, 0) \subset (\C^{d+g+1},0)$ (see \cite{GPRQo}). 
The \textit{local tropicalization} of $(Y,0) \subset (\C^{d+g+1},0)$ may be defined as the closure of the set of tuples 
$(\nu(x_1), \dots, \nu(x_{d+g+1})) \in \R^{d+g+1}$,  where $\nu$ runs through the semivaluations of 
the  ring $\C[[ x_1, \dots, x_{d+1}]]$ such that $\nu (x_i) \in (0, \infty)$, for $i \in \{ 1, \dots, d+g+1 \}$.
Based upon previous results in \cite{MMFGG, GPRQo} and the equivalent definitions of the local tropicalization given in \cite{PS13, PS25}, we
proved that the local tropicalization of $Y$ is equal to the support of a $(d+1)$-dimensional
 fan $\Theta \subset \R^{d+g+1}_{\geq 0}$ determined by the characteristic exponents of $\zeta$.
 The fan $\Theta$ is  isomorphic to the conic polyhedral complex $\bar{\Theta}$ associated with the toroidal embedded normalization $\Psi_g$ (see Corollary  \ref{th-trop}).
This allows us to  reformulate  our main results about the multiplier ideals in terms of the local tropicalization 
of $Y$. We illustrate our results in  Example \ref{Ex:Two pairs qo jumping numbers}.

\medskip

This paper has grown by walking through the beautiful mathematical paths,  full of inspiration and perspectives, which Bernard Teissier has shared with great energy and enthusiasm. We dedicate it to him in the occasion of his 80th birthday.

\section{Basic notions about multiplier ideals and jumping numbers}
\label{Sec: jn}

In this section, we briefly review the basic definitions and properties of the theory of multiplier ideals. For further details, we refer to Lazarsfeld's book \cite[Chapter 9]{PAG}.

\medskip

Let  $Y$ be a  smooth complex
 algebraic variety  and let  $\mathfrak{a}\subseteq \mathcal{O}_{Y}$ be an ideal sheaf. 
 Let 
 $ \Psi :Z\rightarrow Y$
  be a modification, i.e., a proper and birational map,  with $Z$ a smooth variety. Denote by 
$E$ the exceptional divisor of $\Psi$.  The map $\Psi$ is  \textit{log-resolution} 
 of $\mathfrak{a}$ if $\Psi^{*}\mathfrak{a}=\mathcal{O}_{Z}(-F)$, where 
$F$ is an effective divisor such that 
$F+E$ is a \textit{simple normal crossing divisor}. 
 Denote by $\boxed{K_{\Psi}}$ the relative canonical divisor,  that is, the divisor 
associated with the jacobian determinant of $\Psi$.

\medskip

For  $a \in \Q$, we denote by $\lfloor a \rfloor$ the greatest integer lower than or equal to $a$. 
For a $\Q$-divisor  $D = \sum_j {a_j } D_j$,  supported on the prime divisors $D_j$, we denote by $\boxed{\lfloor D \rfloor}\coloneqq \sum_j \lfloor {a_j }  \rfloor D_j$.

\begin{definition} 
\label{Def:Multiplier ideal} \index{multiplier ideals, algebraic}
Let $\Psi :Z\rightarrow Y$ be a log-resolution of an ideal sheaf $\mathfrak{a}$ of $\mathcal{O}_Y$. The \textit{multiplier ideal} sheaf $\mathcal{J}(\mathfrak{a}^{\xi})$
associated to $\xi \in\mathbb{Q}_{\geq 0}$  and $\mathfrak{a}$ is 
$\boxed{\mathcal{J}(\mathfrak{a}^{\xi})} \coloneqq 
	\Psi_{*}\mathcal{O}_{{Z}}(K_{\Psi}-\lfloor  \xi F\rfloor )$.
\end{definition}

The definition of the multiplier ideal $\mathcal{J}(\mathfrak{a}^{\xi})$ relies on the choice of a 
log-resolution of $\mathfrak{a}$, but  it is independent of it (see \cite[Th.  9.2.18]{PAG}).
The multiplier ideal $\mathcal{J}(\mathfrak{a}^{\xi})$ can be expressed  in terms of valuations. If $E_i$ is a prime divisor on $Z$ we denote by 
$\nu_{E_i}$ the vanishing order valuation along $E_i$. 
A prime divisor $E_i$  contained in the support of $F + K_\Psi$ is either the strict transform of a divisor on $Y$ or  it is a component of the exceptional divisor $E$.
Let us write 
$F		=  \sum r_{i}E_{i}$, 
and
\begin{equation} K_{\Psi} = \sum ({\lambda_{E_{i}}}-1) E_{i},  \end{equation}
where the $E_{i}$ are the prime divisors in the support of $E+F$ on $Z$. 
The number $\boxed{\lambda_{{E}_i}}$ is called the \textit{log-discrepancy} of the prime divisor $E_i$. 
It is a birational invariant since it depends only on the associated vanishing order valuation $\nu_{E_i}$.
  Then:
\begin{align} \label{Jac}
\mathcal{J} (\mathfrak{a}^\xi) = 
	\{ h \in \mathcal{O}_{Y} \mid  \nu_{E_i}(h) \geq \lfloor \xi r_{i} \rfloor - (\lambda_{E_{i}}-1) 
		\mbox{ for all } i\}, 
\end{align}
or, equivalently,
\begin{align} \label{eq:form MI}
\mathcal{J} (\mathfrak{a}^ \xi )  = 
	\{ h \in \mathcal{O}_{{Y}} \vert \mbox{ }
	 \nu_{E_i}(h) + \lambda_{E_{i}} > \xi \mbox{ }  r_i  
		\mbox{ for all } i \} .  
\end{align}
The equivalence follows  since for any $a \in \Z$ and $b \in \Q$ one has that $a \geq \lfloor b \rfloor$ if and only if $a > b-1$.

Let $D$ be an effective integral divisor on $Y$. 
Since  $D$ is a Cartier divisor on a smooth variety $Y$, 
 it determines the ideal sheaf   
$\mathcal{O}_Y (-D) =
	\{ h \in \mathcal{O}_Y \mid \mathrm{div} (h) - D \geq 0 \}$. 
The multiplier ideal 	$\mathcal{J} (\mathcal{O}_Y (-D)^\xi )$ 
associated with the ideal sheaf $\mathcal{O}_Y (-D)$ and the number $\xi \in \Q_{\geq 0}$, 
is denoted by $\mathcal{J}(\xi D)$.

 The following lemma introduces the notion of jumping numbers of  the multiplier ideals. 
 \begin{lemma} \cite [Lem. 9.3.21]{PAG}
\label{Def:Jumping numbers}  
Let $Y$ be an smooth algebraic variety and let $\mathfrak{a} \subseteq \mathcal{O}_{Y}$  be a ideal sheaf and 
$o\in Y$.
There exists a strictly increasing discrete sequence $(\xi_{i})$ of positive rational numbers 
such that if  $\xi \in \Q \cap [\xi_{i}, \xi_{i+1})$, then
$
\mathcal{J}(\mathfrak{a}^{\xi_{i}})_o = 
	\mathcal{J}(\mathfrak{a}^{\xi})_o \supsetneq \mathcal{J}(\mathfrak{a}^{\xi_{i+1}})_o . 
$
\end{lemma} 

\medskip 
The numbers $\boxed{\xi_{i}}$ are called the \textit{jumping numbers} associated with $\mathfrak{a}$ at $y$. The smallest jumping number $\xi_{1}$ is called the \textit{log-canonical threshold} of $\mathfrak{a}$ at $y$.
In particular, we have that:

\begin{lemma}  \label{Lem:periodicidad}  \cite [Prop. 9.2.31]{PAG}
Let $D$ be an effective integral divisor on a smooth variety $Y$.  
		Then, 
		$\mathcal{J}(D) = \mathcal{O}_Y (-D)$, and $1$ is a jumping number of the multiplier ideals of $D$.
		In addition, we have that 
		$\mathcal{J} \left( (\xi + 1 ) D \right) = \mathcal{J} \left( \xi  D \right) \otimes \mathcal{O}_Y (-D)$.		
\end{lemma} 
By Lemma \ref{Lem:periodicidad}, which also applies in this local case, 
the jumping numbers of  the multiplier ideals  $\mathcal{J}(\xi D)_y$
 are determined by the  jumping numbers lying in the interval  $(0,1]$
(see   \cite[Ex. 1.7 and Rem. 1.15]{ELSVJumpingCoefficients}).

\medskip

\begin{remark} 
If $Y$ is a smooth variety and $f$ is a germ of complex analytic function at $y \in Y$,
then the above definitions of multiplier ideals and jumping numbers of $(f)$ 
generalize to this local setting (see \cite[Remark 1.26]{ELSVJumpingCoefficients}). 
\end{remark}

\section{Basic notions on toric geometry and Newton polyhedra}     \label{sec-toric}

In this section we introduce the notations and basic results of
toric geometry which we use in this paper (see \cite{FTV,EwCombConv,
TorEmb,OdaCB}).

\medskip

If $u, v \in \Q^d$,  we consider the preorder relation given by
$u \leq v$ if $v \in u + \Q^d_{\geq 0}$. We set also $u < v$
if       $u \leq v$ and $ u \ne  v$. The notation $u \nleq
v$ means that the relation $u \leq v$ does not hold. 
\medskip 

A \textit{lattice} is a free abelian group of finite rank. 
If $\boxed{N}$ is a lattice of rank $d$, i.e., $N\cong \Z^{d}$,  we denote by $\boxed{N_\R}$ the vector
space spanned by $N$ over the field $\R$. 
We denote by $M$ the dual lattice, i.e., $\boxed{M} =\mathrm{Hom} (N,\Z)$
and by 
$
\boxed{\langle \,  , \, \rangle} : N \times M \to \Z 
$
the duality pairing
between the lattices $N$ and $M$.  

\medskip 
In what follows, a {\it
cone} means a {\it rational convex polyhedral cone}: the set $\R_{\geq} v_1 + \cdots +\R_{\geq 0} v_r$ of non
negative linear combinations of vectors $v_1, \dots, v_r \in N$.
The cone $\t$ is {\it strictly convex} if it contains no lines, and  in
such  a case we denote by  $\boxed{0}$ the $0$-dimensional face of $\t$.  The
dimension of $\t$ is the dimension of the linear span of $\t$ in
$N_\R$.  
The cone $\t$ is \textit{simplicial} if the number
of edges of $\t$ is equal to the dimension of $\t$. 
The cone $\t$
is \textit{regular} if it is generated by a sequence of vectors $v_1, \dots, v_r$ 
which is contained in a basis of the lattice $N$.
We denote by $\boxed{N_\tau}$ the sublattice of $N$ spanned by $\tau \cap N$. 
We denote by $\boxed{\mathrm{int}({\tau})}$ the \textit{relative interior} of the cone $\tau$, 
that is, the interior of $\tau$ when viewed on the vector space $N_{\tau,  \R}$.
The {\em dual cone}  $\t^\vee$ and the orthogonal cone $\t^\bot$ of $\t$ are: 
\[
\boxed{\t^\vee} := \{ w  \in M_\R\ |\ \langle u, w
\rangle \geq 0,  \; \forall
u \in \t \}, 
\quad 
\boxed{\t^\bot}:= \{ w  \in M_\R\ | \ \langle u, w
\rangle =  0,  \; \forall
u \in \t \}.
\]

\medskip 

The relation $\theta \leq \t$ (resp. $\theta < \t$)
denotes that $\theta$ is a face of $\t$ (resp. $\theta \ne \t$ is
a face of $\t$).  \medskip

Let ${\s}$  be  a  strictly convex cone rational for the lattice
$N$. By Gordan's Lemma the semigroup ${\s}^\vee \cap M$ is
finitely generated. We denote by $\boxed{\C[{\s}^\vee \cap M]}:= \{ \sum_{\mathrm{finite}}  \ c_\d x^\d \mid
c_\d \in \C, \d \in {\s}^\vee \cap M \}$ the
\textit{semigroup algebra} of $\s^\vee \cap M$ with complex coefficients.
The \textit{affine toric variety} $\boxed{Z_{{\s}} } := \mathrm{Spec} \C [ {\s}^\vee \cap
M]$, denoted also by $\boxed{Z_{{\s}, N}}$, is 
normal. The \textit{torus} $ \boxed{T_N} := Z_{0, N}$ is an open dense subset of
$Z_{\s}$, which acts on $Z_{\s}$ and the action extends the action
of the torus on itself by multiplication. There is a one to one
correspondence between the faces $\t $ of ${\s}$ and the orbits
$\boxed{\orb_\t}$ of the torus action on $Z_{\s}$,
 which reverses the
inclusions of their closures. The orbit $\orb_\t$ is a closed subset of the chart 
$Z_\tau$ and the coordinate ring of  $\orb_\t$ is equal to $\C[\t^\bot \cap M]$. 

\medskip 
If $\theta$ is $d$-dimensional strictly convex cone then its orbit
$\orb_\theta$ is reduced to a closed point called the \textit{origin}
of the toric variety $Z_\theta$. We denote by 
$ \boxed{\C [[\theta^\vee \cap M]] }
=  \{ \sum  c_\d x^\d \mid
c_\d \in \C, \d \in {\s}^\vee \cap M \}$
 the ring of formal power series with exponents in $\theta^\vee
\cap M$. This ring is the completion with respect to the maximal
ideal of the ring $\boxed{\C \{ \theta^\vee \cap M \}} \subset  \C [[\theta^\vee
\cap M]] $ of germs of holomorphic functions at the origin of the toric variety 
$Z_{\theta, N}$.

\medskip

A {\em fan $\boxed{\Sigma}$ of the  lattice $N$} is a family of strictly convex cones  in $N_\R$, which are rational 
for the lattice $N$,
such that for any $\s, \s' \in \Sigma$ we have $\s \cap \s' \in
\Sigma$,  and if $\t \leq \s$ then  $\t \in \Sigma$.  The {\em
support} of the fan $\Sigma$ is the set $\boxed{|\Sigma |} := \bigcup_{\t
\in \Sigma} \t \subset N_\R$. 
 The \textit{$k$-skeletton} of the fan $\Sigma$ is the set $\boxed{\Sigma{(k)}}$ whose elements are the $k$-dimensional cones of $\Sigma$, for $k \in \{0, \dots, d\}$.

\medskip
The
 affine varieties $Z_\s$ corresponding to cones in a fan $\Sigma$
 glue up to define a {\em toric variety} $
 \boxed{Z_\Sigma}$.
If  $\rho$ is a ray of the fan $\Sigma$, we denote  the closure of the orbit  ${\mathrm{orb}_\rho}$  by 
 $\boxed{D_{\rho}} \subset Z_\Sigma$.
 The \textit{boundary} $\boxed{\partial Z_\Sigma}$ of $Z_\Sigma$ is the complement of the torus $T_N$ in $Z_\Sigma$.
Notice that $\partial Z_\Sigma$ 
is equal to the support of the divisor $\sum_{\r \in \Sigma^{(1)} } D_{\rho}$.

 \medskip

A fan $\Sigma'$ is a
{\it subdivision} of the fan $\Sigma$ if both fans have the same
support and if every cone of $\Sigma'$ is contained in a cone of
$\Sigma$. The fan $\Sigma'$ is a \textit{regular subdivision} of
$\Sigma$ if all the cones in $\Sigma'$ are regular and  $\theta
\in \Sigma'$ for any regular cone $\theta \in \Sigma$. Every fan
admits a regular subdivision.
A subdivision $\Sigma'$ of  the fan $\Sigma$ defines a {\em toric
 modification} 
 $
 Z_{\Sigma '}
 \rightarrow   Z_\Sigma
$.
If $\Sigma'$ is a regular subdivision of $\Sigma$, then this
modification 
is a \textit{resolution of singularities} of
$Z_{\Sigma}$.

\medskip 
 
If $u \in |\Sigma|$ we 
denote by  $\boxed{\nu_{u}}$ the \textit{monomial valuation} which satisfies that 
\begin{align} \label{eq:monomial valuation}
\nu_{u} (x^{v}) = 
	\left\langle u, v \right\rangle.
\end{align}
If $u \in N  \cap |\Sigma|$ is a primitive integral vector defining a ray $\rho$
of $\Sigma$,  then $\nu_{u}$ is \textit{vanishing order valuation} along the divisor $D_\rho$.
If $0 \ne u$ does not define a ray of $\Sigma$, we can take a subdivision $\Sigma'$ of $\Sigma$, such that 
the vector $u$ defines a ray 
$\rho \in \Sigma'$, and then $\nu_{u}$ is the \textit{divisorial valuation} associated with
the exceptional divisor $D_\rho$ which appears on $Z_{\Sigma'}$.

\medskip 
Let us fix a strictly convex cone ${\theta}$ of dimension $d$ and rational for the lattice
$N$.
Let $\mathcal{A}  \subset \theta^\vee \cap M$ be a non-empty set. 
The \textit{Newton polyhedron} of $\mathcal{A}$
is the convex hull  $\boxed{\Newton (\mathcal{A})}$  of the set $\mathcal{A} + \theta^\vee$. 
The \textit{support function} $\boxed{\Phi_\mathcal{A}} : \theta \to \R_{\geq 0}$ of $\mathcal{A}$ is
defined by 
$\Phi_\mathcal{A} (u) = \inf
\{    \langle u, v \rangle \mid v \in \mathcal{A}  \}$.

If $\xi >0$ and if $\xi \mathcal{A} := \{ \xi  v \mid v \in \mathcal{A} \}$ then we have 
\begin{equation} \label{eq-sf-hom}
\Phi_{\xi  \mathcal{A}} = \xi \Phi_{\mathcal{A}}. 
\end{equation}
If $w \in \theta^\vee$ and if $w + \mathcal{A}  = \{ w + v \mid v \in \mathcal{A}\}$ then one has that
\begin{equation} \label{eq-sf-tras}
\Phi_{w+ \mathcal{A}} = w + \Phi_{\mathcal{A}}, 
\end{equation}
where  $w \in M_\R$ is seen as linear form on $N_\R$. 

A vector $u \in \theta$  defines the face  $ \boxed{{\mathcal F}_u (\mathcal{A})} :=
\{ v \in {\mathcal N}( \mathcal{A})  \mid \langle u , v
  \rangle   = \Phi_\xi (u  ) \}$
of  the polyhedron ${\mathcal N} (\mathcal{A})$, and all faces of ${\mathcal N} (\mathcal{A})$ are of this form.
We denote by $\boxed{\Sigma (\mathcal{A})}$ the set 
consisting of the cones $ \s( {\mathcal F} ) := \{ u \in \theta
\; \mid \langle u , v \rangle   = \Phi_\xi (u),
 \; \forall v \in {\mathcal F}\}$,
for ${\mathcal F}$ running through the faces of ${\mathcal N}
(\mathcal{A})$.  
The set $\Sigma (\mathcal{A})$ is a  fan of the lattice $N$ supported on $\theta$, which is dual to $\Newton(\mathcal{A})$.

\medskip 

Let us recall how the support function $\Phi_\mathcal{A}$ determines the Newton polyhedron $\Newton (\mathcal{A})$ and vice versa. 
First, notice that 
\begin{equation} \label{eq: New-xi}
\Phi_\mathcal{A} = \Phi_{\Newton (\mathcal{A})}.
\end{equation} 
A vector  $v \in \theta^\vee$ belongs to $\Newton (\mathcal{A})$ if and only if 
\begin{equation} \label{eq: sf-Xi-0}
\Phi_\mathcal{A} (u) \leq \langle u, v \rangle, \mbox{ for all }  u \in \theta.
\end{equation}
In practice, we can check this condition by looking only to the primitive integral vectors $u$ defining rays of the fan 
$\Sigma(\mathcal{A})$. 
This is a reformulation of the description of a polyhedron as the intersection of those half-spaces supporting its facets. 

\begin{proposition} \label{prop: int-pol}
Let $\mathcal{A}$ and $\mathcal{A}'$ be two non-empty subsets of $\theta^\vee \cap M$. 
The following are equivalent: 
\begin{enumerate}[label=(\alph*)]
\item \label{a-xi} The polyhedron $\Newton(\mathcal{A}')$ is contained in the interior of $\Newton (\mathcal{A})$.

\item  \label{b-xi}  For every nonzero vector $u \in \theta$ we have that $\Phi_{\mathcal{A}}  (u)  < \Phi_{\mathcal{A}'}  (u)$.  

\item  \label{c-xi} For every vector $u$ defining a ray of $\Sigma(\mathcal{A})$ we have that  $\Phi_{\mathcal{A}}  (u)  < \Phi_{\mathcal{A}'}  (u)$. 
\end{enumerate}
\end{proposition}
\begin{proof}
It is clear that \ref{b-xi} implies \ref{c-xi}. 

\noindent
 \ref{c-xi} $\Rightarrow$ \ref{a-xi}. 
Take $v \in \mathcal{A}'$.
By hypothesis for all $u \in \theta$ defining a ray of $\Sigma(\mathcal{A})$ 
we have that 
\[
\Phi_{\mathcal{A}}  (u)  < \Phi_{\mathcal{A}'}  (u) =  \inf\{ \langle u , v' \rangle \mid {v' \in \mathcal{A}'} \} \leq \langle u, v \rangle. 
\]
This shows that all inequalities \eqref{eq: sf-Xi-0} hold strictly, thus the vector $v$ belongs to the interior 
of $\Newton (\mathcal{A})$ and \ref{a-xi} holds.

\noindent
 \ref{a-xi} $\Rightarrow$ \ref{b-xi}. Take a nonzero vector $u \in \theta$. 
The inclusion 
$\Newton (\mathcal{A}') \subset \Newton (\mathcal{A})$ implies that 
\[
\Phi_{\mathcal{A}} (u)  \stackrel{\eqref{eq: New-xi}}{=} \Phi_{\Newton (\mathcal{A}) } (u)  \leq   \Phi_{\Newton (\mathcal{A}') } (u)  \stackrel{\eqref{eq: New-xi}}{=}   \Phi_{\mathcal{A}'} (u). 
\]
Equality holds in this formula precisely when the face ${\mathcal F}_u (\mathcal{A}')$ of $\Newton (\mathcal{A}')$  is contained in  ${\mathcal F}_u (\mathcal{A})$. 
Since $u \ne 0$, the face  ${\mathcal F}_u (\mathcal{A})$  is contained in the boundary of  $\Newton (\mathcal{A})$. 
By hypothesis,  we get that ${\mathcal F}_u (\mathcal{A}') \cap {\mathcal F}_u (\mathcal{A})   = \emptyset$, hence 
the inequality is strict.
This shows that \ref{b-xi} holds. 
\end{proof}

 \begin{notation}\label{newton}
The support of a a nonzero series $h= \sum
c_\d x^{\d} \in \C \{ \theta^\vee \cap M \}$
is the set $\boxed{\mathrm{supp} (h)} = \{ \d \in M \mid c_\d \ne 0 \}$. 
We write also $\boxed{\Phi_h}, \boxed{\Newton(h)} $,  $\boxed{\Sigma(h)}$ or $\boxed{{\mathcal F}_u (h)}$ instead of 
$\Phi_\mathcal{A}, \Newton (\mathcal{A})$, $\Sigma(\mathcal{A})$  or ${\mathcal F}_u (\mathcal{A})$ respectively,
when $\mathcal{A} = \mathrm{supp} (h)$. 
\end{notation}

\begin{remark}
Notice that if  $u \in \theta$ we have that 
\begin{equation} \label{eq: Phi-h-0}
\nu_{u} ( h ) = 
	\Phi_{h} (u).
\end{equation}
\end{remark}

In the following definition we follow the terminology of \cite[Section 1.4.1]{GBGPPP20}.
\begin{definition} \label{def-Newtonm}
Given a nonzero series  $h= \sum
c_\d x^{\d} \in \C \{ \theta^\vee \cap M \}$ we say that: 
\begin{itemize}
\item  The fan $\Sigma(h)$ is the \textit{Newton fan} of $h$.
\item The toric modification 
$
 \boxed{\psi(h)}:  Z_{\Sigma(h)} \longrightarrow Z_{\theta}  
$
defined by the fan $\Sigma(h)$ is 
\textit{Newton modification} determined by $h$.
\end{itemize}
 \end{definition}
 
\begin{lemma} \label{lem-0-dim}
The strict transform of the germ defined by $h$  by the Newton modification 
$\psi(h)$ does not pass through the $0$-dimensional orbits  of the toric variety $Z_{\Sigma(h)}$. 
\end{lemma}
\begin{proof}
Denote by $H$ the germ defined by $h =0$.
A $0$-dimensional orbit is the origin of a chart $Z_{\theta}$ where
$\theta \in \Sigma(h)$ is the dual face of a vertex $m$ of $\mathcal{N} (h)$. 
This implies that on $Z_{\theta}$ we can factor 
$
h = x^{m} (c_m + \dots )$, 
where  $h x^{-m}$ defines the strict transform of $H$ on $Z_\theta$. 
Since $m$ is a vertex of  $\mathcal{N} (h)$, we have that $c_m \ne 0$, hence 
the strict transform of $H$ does not pass through the origin of $Z_\theta$.
\end{proof}

\subsection{Toroidal embeddings} \label{sec-TE}
The theory of toroidal embeddings
is a generalization of the theory of toric varieties. 
We outline some definitions  and  results and we refer to
the classical book of Kempf \textit{et al.~}\cite{TorEmb} for details.

\medskip

Let $V$ a normal complex variety of dimension $d$ and $U \subset V$ a smooth Zariski open subset of $V$. 
We say that $\boxed{U \subset V}$ is a  \textit{toroidal embedding} if for every closed point $v \in V$ there exists 
a $d$-dimensional affine toric variety $Z_{\s_v, N}$ of dimension $d$, a closed point $z_v \in Z_{\s_v, N}$, and an 
isomorphism of complete local $\C$-algebras: 
\[
\hat{\mathcal{O}}_{V,v} \longrightarrow \hat{\mathcal{O}}_{Z_{\s_v, N}, z_v},
\]
such that the ideal of $\hat{\mathcal{O}}_{V,v}$ generated by the ideal of $V \setminus U$ corresponds under 
this isomorphism to the ideal of $ \hat{\mathcal{O}}_{Z_{\s_v, N}, z_v}$, generated by the ideal of $Z_{\s, N} \setminus T_N$. 
By restricting $Z_{\s_v, N}$ to an open equivariant subset containg $z_v$, we assume that 
the orbit of $z_v$ is closed in $Z_{\s_v, N}$.
This isomorphism is called a  \textit{local model} at $v$.
\medskip

If $U \subset V$ is a toroidal embedding then we have a \textit{boundary divisor} 
$ \boxed{\partial V}: = V - U$, which decomposes as $\partial V = \bigcup_{i \in I} D_i$,   where the prime components  $D_i$ are of codimension one. 
The toroidal embedding $U \subset V$ is \textit{without self-intersection} if $D_i$ is normal for all $i \in I$.
We will assume this hypothesis in what follows. 

\medskip
The toroidal embedding without self-intersection $U \subset V$ defines a stratification of $V$ with strata given by the irreducible components of the sets $\bigcap_{i \in J} D_i \setminus \bigcup_{i \notin J} D_i$, for $J \subset I$. If $Y$  is a
stratum of this stratification the star of $Y$, denoted by $\boxed{\mathrm{Star} \, Y}$,  is the open subset defined as the union  of the strata $Z$ such that $Y \subset \overline{Z}$.

\medskip 

If $Y$ is a stratum of this stratification we denote by $\boxed{M^Y}$ the group of Cartier divisors on $\mathrm{Star} \, Y$, supported on 
$\mathrm{Star} \, Y \setminus U$ and by $\boxed{N^Y}$ its dual group. We have that 
$N^Y$ and $M^Y$ are free abelian groups of rank equal to the codimension of $Y$.
We denote by $\boxed{M^Y_+}$ the subsemigroup of $M^Y$ consisting of effective divisors. Then, the cone 

\[
\boxed{\s^Y} := \{ u \in N^Y_\R \mid \langle D, u \rangle \geq 0 , \, \forall D \in M_+^Y \} \subset N^Y_\R,
\]
 is rational with respect to the lattice $N^Y$, polyhedral and of dimension equal to the rank of $N^Y$.
If $v \in Y$,  the cone $\s_{v}$ of a local model at $v$ and the sublatice  $N_{\s_v}$ spanned by 
the vectors in $N \cap \s_{v}$,  are independent of the choice of $v \in Y$.  Namely, we have canonical identifications $\s_v \cong \s^Y$  and $N^Y \cong N_{\s_v}$.

\medskip 

Let $U \subset V$ be a toroidal embedding without self intersection. 
It has an associated  \textit{conic polyhedral 
complex with integral structure} $\boxed{\Xi}$, consisting of the cones $\s^Y$ of the lattice $N^Y$, for $Y$ running through the strata. The support  $\boxed{|\Xi|}$ of the complex $\boxed{\Xi}$ is a topological space 
obtained from the disjoint union of the cones $\s^Y$ under a natural equivalence relation. Namely, 
if $Z \subset \mathrm{Star} \, Y $ is a stratum of the stratification,  one has an associated 
face $\tau^Z$ of $\s^Y$, and natural identifications of the pairs $(\sigma^Z, N^Z)$ and $(\tau^Z,  N^Y_{\tau^Z})$.
In particular, we have a one to one correspondence which sends a ray $\theta$ of
$\Xi$ to  the prime component  $\boxed{D_\theta}$  of $\partial V$, which is equal to the closure of the correspondent stratum.

\medskip 
A  finite rational polyhedral subdivision
$\Xi'$  of $\Xi$ induces a modification $V' \rightarrow V$. This modification 
defines an isomorphism $U' \to U$ over $U$, for some open subset $U' \subset V'$,
in such a way that  the pair  $(V', U')$ is also a toroidal embedding. In particular, if  $\Xi'$ is a \textit{regular
subdivision} of $\Xi$
then the map $V' \to V$ is a resolution of singularities of $V$.

\medskip 

For instance, if $V = Z_{\Sigma}$ is the toric variety of the fan $\Sigma$ of the lattice $N$ then
$T_N \subset Z_\Sigma$
a toroidal embedding without self-intersection. Its associated conic polyhedral complex is determined by the fan $\Sigma$.   In comparison with the toric varieties, the cones of a fan are contained inside the same vector space, but in general if $Y$ is an stratum of the toroidal embedding $U \subset V$,  the cone $\s^Y$ is a subset of the vector space $N^Y_\R$, which depends on the stratum $Y$.

\section{Quasi-ordinary hypersurface singularities}\label{qoh}

In this section we recall the definition and the main properties of quasi-ordinary hypersurface singularities following \cite{A55, LipmanQOSingEmbSurf, LipmanTopInvQO, Gau, PPPDuke, SgrPG, GPRQo}. 

\medskip

A germ $\boxed{(\X, o)} $ of  equidimensional complex analytic variety  of dimension $d$
 is \textit{quasi-ordinary} (q.o.) if there exists a
finite projection $\boxed{\p}: (\X, o) \rightarrow (\C^d,0)$ which is a
local isomorphism outside a germ of normal crossing divisor.
In the hypersurface case  $(\X, o) \subset (\C^{d+1}, 0)$, there exists  local coordinates  
$(x_1, \dots, x_d, y)$ at the origin of $\C^{d+1}$, such that 
$\pi$ is the restriction to $(\X,o)$ of
the projection 
\begin{equation}
\label{projection} (\C^{d+1},0)  \to (\C^d,0),  \quad (x_1, \dots, x_d,
y) \mapsto (x_1,\dots, x_d).
\end{equation}

\medskip
Let us set $ \x  = (x_1, \dots, x_d)$. 
If $h$ belongs to the ring of convergent power series
$\C \{ \x , y \}$ we denote by $\boxed{Z(h)}$
the germ defined by $h= 0$ at the origin of $\C^{d+1}$.

\medskip 

{A reduced power  series $\boxed{F} \in \C \{ \x, y \}$,  is called \textit{quasi-ordinary} if the germ $Z(F)$ is quasi-ordinary with respect to the projection \eqref{projection}.
This implies that the order of $F({0},y)$ is finite.
By the Weierstrass Preparation Theorem 
there is a unit $U \in  \C \{ \x, y \}$ such 
that $\boxed{f} := U\cdot F$ is a Weierstrass polynomial in the ring
$\C\{ \x  \} [{y}]$. We say that $f$ is a  \textit{q.o.~polynomial} defining 
the germ $(\X,o)$. 
The hypothesis implies that the discriminant 
$\boxed{\D_{{y}} f}$ of $f$  is the product of a monomial in $\x$ and 
a unit in 
$\C \{ \x \}$.

 \medskip 
 
 In addition, we assume in this paper that the germ $(\X,o)$ is an analytically irreducible, that is,
the polynomial  $f $ is irreducible in  
$\C \{ \x 
\} [{y}]$.
The Abhyankar-Jung theorem asserts that the roots of the q.o.~
polynomial $f$, called {\it q.o.~ branches}, are fractional power
series in the ring $\C \{ \x^{1/n} 
\}$, 
where $\boxed{\x^{1/n}} = (x_1^{1/n}, \dots, x_n^{1/n})$ and 
$n$ is  the degree $\boxed{\deg_y (f)}$  of $f$ in $y$ (see
\cite{A55}, \cite{PR12}).
We denote the ring of  generalized \textit{Puiseux series} in the variables $\x^{1/n}$ by 
\[
\boxed{\C \{ \x^{1/ \N}\}} := \bigcup_{m \geq 1 } \C \{ \x^{1/m} \}.
\]

The roots of quasi-ordinary polynomials are special among the generalized Puiseux series. 
If $\z_i$ and $\z_j$  are two different roots of the quasi-ordinary polynomial $f$ in the ring 
$\C \{  \x^{1/n} \}$
then 
$\z_i - \z_j$ divides the 
discriminant $\D_{{y}} f$. 
This implies that 
$
\z_i - \z_j = \x^{ \a_{i, j} } \cdot u_{i, j},
$
where $ \a_{i, j} \in \frac{1}{n} \Z^d_{\geq 0}$, and $u_{i, j}$ is a unit in $\C \{  \x^{1/n}  \}$. 
The fractional monomials $x^{\a_{i,j}}$ are called the \emph{characteristic 
monomials} of the  q.o. polynomial $f$ (or the q.o. series $F$). 
The tuples $\a_{i, j }$, for $i \ne j$,  are called the generalized \emph{characteristic 
exponents} of $f$. These tuples can be relabelled as $\boxed{\a_1, \dots, \a_g}$
in such a way that: 
\begin{equation} \label{eq:ch-exp}
 0 <  \,  \a_1 < \a_2 < \cdots < \a_g,
\end{equation}
where $<$ means $\ne$ and $\leq$ coordinate-wise. See \cite{LipmanQOSing, LipmanTopInvQO}.

\begin{definition} \label{char-lattices}
We consider the sequence of \emph{characteristic lattices}
\[
\boxed{M_0} := \Z^d, \quad \boxed{M_j} :=
M_{j-1} + \Z \a_j  \mbox{, for } j=1, \dots, g,
\]
and \emph{characteristic indices}
\[
\boxed{n_j} := [M_{j-1} : M_j], \mbox{ for } j=1, \dots, g.
\]
We set also 
\[
\boxed{n_0} :=1, \quad \boxed{e_0} := n_1 \cdots n_g   \quad \mbox{ and } 
\quad \boxed{e_{j}} := e_0 / (n_1 \dots
n_j) \mbox{ for } j=1, \dots, g,
\]
and
\begin{equation}\label{rel-semi}
\boxed{{\g}_1} :=  \a_1 \quad \mbox{ and } \quad \boxed{{\g}_{j+1} } :=  n_j {\g}_{j}
+ \a_{j+1} -  \a_{j}, \quad \mbox{ for } \quad j= 1, \dots, g-1,
\end{equation}
We consider also the semigroup:
\[ \boxed{\Gamma}  := \Z^d_{\geq 0} +   \Z_{\geq 0}  \g_1
+ \cdots + \Z_{\geq 0}  \g_g \subset \Q^d_{\geq 0}. \]
By convenience we denote $\boxed{\a_0 } = \boxed{  \g_0 } = 0 \in M$. 
For $\ell \in \{0, \dots, g\}$ we set 
\[\boxed{\Gamma_\ell}
:= \Z^d_{\geq 0} +   \Z_{\geq 0} \g_1
 + \cdots + \Z_{\geq 0} \g_\ell.   \]
\end{definition}
\begin{remark} 
We have that  $n_j \geq 2$, for $j=1, \dots, g$, and $e_0 = n= n_1 \dots n_g$  (see \cite{LipmanTopInvQO}). 
\end{remark}

\medskip 
Let us fix a root  $ \zeta= \zeta(x^{1/n})$ of the q.o. polynomial $f$. 
Denote by $\mathcal{C}(\zeta)$  the set of functions $h \in \C\{ x, y\}$ \textit{with a dominating exponent $\gamma(h)$
 with respect to $\zeta$}, that is, 
$h(x, \zeta(x^{1/n})) = x^{\gamma (h)} \cdot u(h)$, where $\g(h) \in \frac{1}{n} \Z_{\geq 0}^d$ and 
$u(h) \in \C \{ x^{1/n} \}$ is a unit. 
\begin{theorem} \cite{SgrPG, PPPDuke, KiMiSgr}
We have that 
$
\Gamma = \{ \gamma(h) \mid h \in\mathcal{C}(\zeta))\} $. 
\end{theorem}
The  semigroup  $\Gamma$ is an analytic invariant of $D$, which is called the \textit{semigroup of the quasi-ordinary hypersurface} (see  \cite{SgrPG, PPPDuke, KiMiSgr}).
By the results of Gau and Lipman \cite{LipmanTopInvQO,Gau}, the semigroup $\Gamma$  classifies the embedded topological type of the germ $( D , 0) \subset (\C^{d+1}, 0)$, that is, it is a topological invariant.

\medskip

The following lemma provides   distinguished way of writing the elements of the semigroup $\Gamma$:
\begin{lemma} \label{lem:unique-sg} 
\cite[Lem. 3.3]{SgrPG} Take $\ell \in \{ 1, \dots, g \}$.
\begin{enumerate}
\item We have that
$n_\ell \g_\ell \in \Gamma_{\ell-1}$. 
\item Any element in the semigroup $\g \in \Gamma_\ell$ can be represented in a unique way 
in the form 
$
\g = u + i_1 \g_1 + \cdots + i_{\ell} \g_\ell$,
where $u \in \N^d$ and $0 \leq i_j < n_i$ for $j \in \{ 1, \dots, \ell \}$. 
\end{enumerate}
\end{lemma} 
\begin{definition} \label{def:semi-root}
A \textit{semi-root of $f$ of degree} $n_0 \dots n_{j-1}$, for $j \in \{1, \dots, g\}$, is a 
quasi-ordinary power series 
 $ h \in \mathcal{C}(\zeta) $ such that the order of $h (0, y)$ is equal to  $n_0 \dots n_{j-1}$, and
  $\g(h) = \g_{j}$.
We say also that $f$ is the semi-root of $f$ of degree $n = n_0 \dots n_g $. 
\end{definition}

Take  $j \in \{ 0, \dots, g-1\}$.
A semi-root of $f$ of degree $n_0 \dots n_{j}$
can be obtained in terms of the expansion of a 
fixed root $\boxed{\z} \in \C \{ \x^{1/n} 
\}$ of $f$ as follows. We denote by $\boxed{\tau_j}$  the series obtained by eliminating all the terms $x^\alpha$ which have exponent $\alpha \geq \a_{j+1}$ in the expansion of $\z$.
The minimal polynomial $\boxed{g_j}$ of $\tau_j$ over the field of fractions of $\C \{ \x \}$ 
has degree $n_0 \dots n_{j}$, 
belongs to $\C \{ \x \} [{y}]$ and  is a semi-root of $f$  (see \cite{SgrPG}). 

\medskip

In the rest of this paper we fix a semi-root  $\boxed{\x_{d+\ell+1}} \in \C \{ \x \} [{y}]$ of 
degree $n_0 \cdots n_\ell$ of the q.o. polynomial $f$, for  $\ell \in \{ 0, \dots,g \}$. 
\begin{definition}  \label{def: semi}
We say that $ \x_{d+1}, \dots, \x_{d+g+1}$ is the sequence of \emph{characteristic semi-roots} of $f$
and that 
${\x_1, \dots, x_d, \x_{d+1}, \dots, \x_{d+g+1}}$ is a \emph{complete sequence of 
semi-roots} of $f$.
We denote 
the hypersurface germ defined
by $\x_j = 0$ by 
\begin{equation} \label{eq-Dj}
\boxed{\X_{j}} := Z(x_j), \mbox{ for } j =1, \dots, d+ g +1.
\end{equation}
\end{definition}

By \cite {SgrPG} we have that $x_{d +\ell +1}$ is a quasi-ordinary series with 
characteristic exponents $\a_1, \dots, \a_{\ell }$, for $\ell \in \{ 1, \dots,g \}$. 
This implies that the semigroup of of the q.o. series $x_{d +\ell +1}$ is equal to 
 $\Gamma_\ell$. 

\begin{remark} \label{rem:y} 
By Definitions  \ref{def: semi} and \ref{def:semi-root} we have that $x_{d+1} (0, y)$ is of order one in $y$, and then
$(x_1, \dots, x_{d+1})$ is a system of local coordinates at the origin of $\C^{d+1}$.
 \end{remark}

\subsection{Expansions in term of the semi-roots}
In this section, we describe how to expand a series $h \in \C \{\x ,x_{d+1} \}$ 
in terms of the fixed elements $x_1, \dots, x_{d+g+1} \in \C\{ \x \}[x_{d+1}] $ (see Definition \ref{def: semi}).

\medskip

Let $F \in  \C \{  \x \}[y]$  be a  Weierstrass polynomial in $y$.
Let us recall how to build the $F$-\emph{adic expansion} of  a series in $\C \{\x,y\}$.

\begin{lemma} \label{Lema:Weierstrass division}
Let 
$F \in  \C \{  \x \}[y]$ be a  Weierstrass polynomial  of degree $m \geq 1$ in $y$.
Then, any $h \in \C \{\x,y\}$ has a unique $F$-adic expansion of the form:
\begin{equation} \label{f-adic}
h = \sum_{j \geq 0} P_j F^{j} ,  
\end{equation}
such that for any $j \geq 0 $,  the coefficient  $P_j \in \C \{\x \}[y]$ is a polynomial  of degree $<m$.
\end{lemma}

\begin{proof}
The assertion is clear if $h =0$, so assume that $h \ne 0$. 
Since the ring $\C \{\x,y\}$ is a unique factorization domain, 
there exists an integer $s\geq 0$ such that 
\begin{align} \label{f-adic1}
h = h_{s-1} \cdot  F^{s} , \quad \mbox{ with } h_{s-1} \notin (F). 
\end{align}
We apply the Weierstrass division theorem 
(see \cite[Chapter VII, Theorem 5]{ZarSamComAlgVol2} and \cite[Theorem 2.2]{Wolfgang}). 
There are are unique elements $h_{s} \in \C \{\x,y  \}$
and $P_s \in \C \{\x \}[y]$ of degree $<m$ 
such that 
$
h_{s-1}  = h_{s} \cdot F+ P_s$. 
Substituting in (\ref{f-adic1}) we obtain
$
h = P_s \cdot F^s +  h_s \cdot F^{s+1}$.
By applying the Weierstrass division theorem to the quotient $h_s$ 
and by iterating this process, we get for $q > s$ that 
\begin{align*}
h = P_s \cdot F^s + P_{s+1} \cdot F^{s+1} + \cdots + P_{q} \cdot F^{q} +  h_{q+1} \cdot F^{q+1}, 
\end{align*}
where 
$P_j \in \C \{\x \}[y]$ of degree $<m$ for $j = s, \dots, q$ and $h_{q+1} \in  \C \{\x,y  \}$.
Since the order of $F$ is $\geq 1$, setting $P_j =0$, for $j <s$, we have that the equality \eqref{f-adic} holds in the ring $\C [[ \x , y ]] $. The unicity follows from that of the 
Weierstrass division theorem.
\end{proof}

In the following lemma we consider expansions in terms of 
the complete sequence of semi-roots 
$\x_1, \dots, \x_d, \x_{d+1}, \dots, \x_{d+g} \in \C \{ \x \} [x_{d+1}]$ of the q.o. polynomial $f$ (see Definition \ref{def: semi}).

\begin{lemma} {\cite[Lemma 7.2]{PPPDuke}}  \label{Lem: semi-1} 
Let 
$P \in \C \{ \x \} [x_{d+1}]$ be a polynomial of degree $< n$. 
Then, we have a unique expansion of the form: 
\begin{align} \label{Form:Expansionsemi-roots}
P =\sum_I  c_{I}  \, x_1^{i_1}  \cdots x_{d+g}^{i_{d+ g}}, \mbox{ with } c_I  \in \C,
\end{align}
such that if   $c_{I} \ne 0 $ for $I = (i_1, \dots, i_{d+g})$  then
\begin{align} \label{Form:Exponentssemi-roots}
0\leq i_{d+j}< n_{j}\ \mathrm{for}\ j \in \{ 1,\ldots ,g \}.
\end{align}
\end{lemma}

By combining the previous lemmas  we obtain: 
\begin{lemma} \label{Lema:Expansion semi-roots}
Any power series $h \in \C \{\x ,x_{d+1} \}$ 
has a unique expansion: 
\begin{align} \label{Form:Expansion of any germ in terms of semi-roots}
h = {\sum}_{I} 
c_{I} \, 
x_1^{i_1}  \cdots   \x_{d+g+1}^{i_{d+g+1}} , \mbox{ with } c_I \in \C,
\end{align}
such that if $c_I \ne0$, for 
$I =(i_1, \dots, i_{g+1}) $
then
\eqref{Form:Exponentssemi-roots} hold.
In addition, the degrees as polynomials in $x_{d+1}$ of the terms
$\x_{d+1}^{i_{d+1}} \cdots  
 \x_{d+g+1}^{i_{d+g+1}}$ 
are pairwise distinct. 
\end{lemma}

\begin{proof}
We apply first  Lemma \ref{Lema:Weierstrass division} to $h$ with respect to the polynomial $\x_{d+g+1}= f$. 
We obtain an expansion 
$
h = \sum_{m \geq 0} P_m \, \x_{d+g+1}^m$, 
where $P_m \in  \C \{\x \}[x_{d+1}]$, \mbox{ and }  $\deg_{\x_{d+1}} P_m < n$ for every $m \geq 0$. 
The result follows by expanding the coefficients $P_m$ by using Lemma \ref{Lem: semi-1}. 
The assertion about the degrees in $\x_{d+1}$ of the terms of the expansion is elementary
 (see the proof of \cite[Lemma 7.2]{PPPDuke}). 
\end{proof}

\begin{definition} \label{def: expansion}
The expansion \eqref{Form:Expansion of any germ in terms of semi-roots} 
is the \textit{$(x_1, \dots, x_{d+g+1})$-expansion} of $h \in \C \{\x ,x_{d+1} \}$. 
\end{definition}

Notice that the terms  $x_1^{i_1} \cdots x_{d+g+1}^{i_{d+g+1}}$ 
of the expansion \eqref{Form:Expansion of any germ in terms of semi-roots} 
are  \textit{generalized monomial in the semi-roots}. These terms are not monomials in the power series ring
 $\C\{ x_1, \dots, x_{d+1}  \}$ in general.

\subsection{Quasi-ordinary hypersurfaces relative to a toric base}

We explain a generalization of quasi-ordinary hypersurfaces singularities which arise when we 
replace the target $\C^d$ of a quasi-ordinary projection by an affine toric variety. 
This class of singularities was introduced in \cite{PGSQOT}.  It was applied in
 \cite{GPRQo} to describe the toroidal embedded resolutions of quasi-ordinary hypersurfaces  
 (see Section    \ref{resolution}).

\medskip

Let $N'$ be  a lattice of  rank $d+1$, endowed with a direct sum decomposition
$N'= N \oplus \Z \epsilon$, where $\epsilon \in N'$ is a primitive vector.
Then, if $M'$ and $M$ are the dual lattices of $N'$ and $N$ respectively, we also have  a direct sum decomposition 
$M' = M \oplus \Z \varepsilon$, where $\varepsilon$ is the primitive integral vector orthogonal to $N$ such that 
$\langle \epsilon, \varepsilon \rangle =1$.  Let $\r \subset N_\R$ be a strictly convex rational cone of dimension $d$ and
set $\varrho = \r + \R_{\geq 0} \epsilon$.  Then, the dual cone of $\r$ seen in $N_\R$,  is equal to the face of 
$\varrho^\vee \cap \epsilon^\perp$ defined by the vector $\epsilon$. 
We get the inclusion of rings
$
\C[\r^\vee \cap M]  \subset \C[\varrho^\vee \cap M']
$
which corresponds to the toric projection 
\begin{equation} \label{def: pibase}
 Z_{\varrho, N'} \longrightarrow Z_{\r, N}.
\end{equation}

We have that the ring of germs of holomorphic functions at the origin of $Z_{\varrho, N'}$ is equal to 
\begin{equation} \label{eq:rep-local}
\mathcal{O}_{Z_{\varrho, N'}, 0} =  \C\{ \varrho^\vee \cap M' \} = \boxed{ \C\{ \r^\vee \cap M\}\{x^\varepsilon\} }.
\end{equation}

Let $H \in \C\{ \varrho^\vee \cap M' \}$ define the germ $Y = Z(H)$ at the origin of  $Z_{\varrho, N'}$. 
We say that the germ $Y$, or the series $H$,  is \textit{quasi-ordinary relative to the base}  $Z_{\r, N}$  if the
restriction $\pi$ of the projection \eqref{def: pibase} to some representative of the germ $Y$ is unramified over the 
torus $T_N$ of $Z_{\r, N}$.  The restriction  of $H$ to the closure of the orbit defined by $\rho$ in 
$Z_{\varrho, N'}$, is a power series in the ring $\C\{ x^{\varepsilon} \}$. The order $m$ of this series
is the degree of the branched covering defined by $\pi$. 
One has a form of the Weierstrass Preparation Theorem, which guarantees that 
there exists a unit $U \in \C\{ \varrho^\vee \cap M' \}$ and a polynomial $h$  of degree $m$ in the ring 
$\C\{ \r^\vee \cap M\}[x^{\varepsilon}]$ such that $U\cdot H = h$. 
The discriminant of $h$ with respect to $x^{\varepsilon}$ is a monomial times a unit in the ring 
$\C\{ \r^\vee \cap M\}$. If the toric variety $Z_{\r, N}$ is smooth we are in the case of  quasi-ordinary hypersurface singularities considered in the previous section.

\begin{remark}
There is a version of the Abhyankar-Jung Theorem for this class of singularities  (see \cite{PGSQOT} for a topological proof or \cite{PR12} for an algebraic proof). This version states that there exists  a finite lattice extension $M \subset \tilde{M}$ such that the polynomial $h$ has all its roots in the ring $\C \{ \r^\vee \cap \tilde{M} \}$. These roots are called toric quasi-ordinary branches, relative to the toric base $Z_{\r, N}$. The notion of characteristic exponents of a quasi-ordinary branch relative to the base $Z_{\r, N}$ extends with the same definition.
If $\z_i$ and $\z_j$ are two different roots of $h$ in the ring $\C \{ \r^\vee \cap \tilde{M} \}$ then 
the difference $\z_i - \z_j$ divides the discriminant. This implies that there exists a monomial $\a_{i, j} \in  \r^\vee \cap \tilde{M}$ and a unit $u_{i,j} \in \C \{ \r^\vee \cap \tilde{M} \}$ such that:
$
\z_i - \z_j = x^{\a_{i, j}} \cdot u_{i,j}$.
 The \textit{characteristic exponents} $\a_{i,j}$ obtained in this way 
 are vectors in the vector space $M_\Q$, which belong to the \textit{reference
cone} $\r^\vee$.  If the polynomial $h$ is irreducible,  then the characteristic exponents
can be relabeled as $\a_1, \dots, \a_r$ in such a way that 
\[
\a_{i+1} \in \a_i + \rho^\vee, \mbox{ for } i =1, \dots, r-1. 
\]
This relation is a generalization of  \eqref{eq:ch-exp} in the classical case. The notions of characteristic exponents, semigroup and characteristic semi-roots
generalize to the case of of q.o.~ hypersurface singularities relative to a toric base (see \cite{GPRQo} for details).

\end{remark}

\subsection{Some combinatorial facts}
We give in this section some combinatorial definitions and  results in terms of the characteristic exponents 
$\a_1, \dots, \a_g$ of 
a quasi-ordinary hypersurface.

\begin{notation} \label{not:basic}
We denote by 
$\boxed{\epsilon_1, \dots, \epsilon_{d+1}}$  the canonical basis of  
$N_0' : = \Z^{d+1}$, and by  $\boxed{\varepsilon_1, \dots, \varepsilon_{d+1}}$
 the dual basis of the dual lattice $M_0'$.
We denote by $N_0$ the lattice with basis $\epsilon_1, \dots, \epsilon_{d}$,  and by $M_0$ its dual lattice. 
We also define  the cones
\[
\boxed{\varrho} : = \R_{\geq 0} \epsilon_1 + \cdots +\R_{\geq 0} \epsilon_{d+1} \mbox{ and } 
\boxed{\r} = \R_{{\geq 0}} \epsilon_1 + \cdots +\R_{\geq 0} \epsilon_{d}.
\] 
 \end{notation}
 The dual cone of $\varrho$ is 
 $
 \varrho^\vee = {\R_{\geq 0}} \varepsilon_1 + \cdots +\R_{\geq 0} \varepsilon_{d+1}.
 $
The dual cone of $\r$, seen as a cone in $N_\R$, is 
$
\r^\vee = {\R_{\geq 0}} \varepsilon_1 + \cdots +\R_{\geq 0} \varepsilon_{d}.
$
We can view $\r^\vee$ as a face of $\varrho^\vee$. 

\medskip 
We consider the lattices $\boxed{N_\ell'} := N_\ell  \oplus \Z \epsilon_{d+1}$ and $\boxed{M_\ell'} := M_\ell  \oplus \Z \varepsilon_{d+1}$.  Recall that $M_\ell$ was introduced in Definition \ref{char-lattices} in terms of 
the characteristic exponents $\a_1, \dots, \a_g$. Then, the lattice
\begin{equation} \label{eq:Nele}
\boxed{N_\ell} := \{ \nu \in N_0  \mid \langle \nu, \a_j \rangle \in \Z, \mbox{ for } j=1, \dots, \ell \} 
\end{equation}
is the dual lattice of $M_\ell$, for $\ell \in \{0, \dots, g\}$. 
Notice that the semigroup  $ \r^\vee  \cap M_\ell$ may not be generated by 
a basis of the lattice $M_\ell$, for 
$\ell \in \{1, \dots, g\}$.

\begin{definition} \label{def:SigmaL}
For $\ell \in \{1, \dots, g\}$ we introduce the following cones: 
\begin{equation} \label{eq:Sigma-ell}
\boxed{\bar{\s}_\ell^+}
= 
 \{ \nu + r \epsilon_{d+1}  \in \varrho \mid  \langle \nu, \a_\ell - \a_{\ell-1}  \rangle \leq r \}  \mbox{ and }
\boxed{\bar{\s}_\ell^-}
= 
 \{ \nu + r \epsilon_{d+1}  \in \varrho \mid   r    \leq  \langle \nu, \a_\ell- \a_{\ell-1}  \rangle \}, 
\end{equation}
which intersect along the $d$-dimensional cone 
\[
\begin{array}{ccl}
\boxed{\bar{\r_\ell}}  & =  &  \{ (\nu, r) \in \varrho \mid   r   =  \langle \nu, \a_\ell - \a_{\ell -1} \rangle \}. 
\end{array}
\]
We have a subdivision  $\boxed{\Sigma_\ell}$ of the cone $\varrho$ such that its $(d+1)$-dimensional skeleton is 
$\Sigma_\ell (d+1) = \{ \bar{\s}_\ell^+ , \bar{\s}_\ell^- \}$.
We define $\boxed{\Sigma_{g+1}}$ as the fan of faces of the cone $\varrho$ with respect to the lattice $N_g'$.
\end{definition}

In Notation \ref{cones} below we will introduce rational cones $\s_\ell^+, \s_\ell^+, \rho_\ell \subset \R^{d+g+1}_{\geq 0}$, for $\ell \in \{1, \dots, g\}$, and then we will explain the relations between these cones 
and those introduced in Definition \ref{def:SigmaL}. We follow here the notation introduced in \cite{MMFGG}.

\begin{lemma} For $\ell \in \{1, \dots, g \}$, there is an isomorphism of toric varieties 
\begin{align} \label{eq:rhoj-3}
Z_{{\rho_{\ell}},N_{\ell-1}'} \longrightarrow 
	Z_{\rho , N_\ell} \times \orb_{\rho_\ell,N_{\ell-1}'}  = Z_{\rho , N_\ell} \times \C^*.
\end{align}
\end{lemma}
\begin{proof}

For $\ell \in \{ 1, \dots, g\}$ we consider the  lattice homomorphism 
\begin{equation} \label{eq:phi-ell}
\boxed{\phi_\ell}: M'_{\ell-1} \to M_\ell, \quad 
\sum_{i = 1}^{d+1} m_i \varepsilon_i \mapsto    \sum_{i = 1}^{d} m_i \varepsilon_i   + m_{d+1} (\a_\ell - \a_{\ell-1}).
\end{equation}
By Definition \ref{char-lattices},  the homomorphism $\phi_\ell$ 
is surjective. It has kernel
\[\ker
(\phi_\ell) =   \left(
n_\ell \varepsilon_{d+1} -n_\ell ( \a_\ell- \a_{\ell-1}) )  \right) \Z,\] 
since $n_\ell$  is equal to the index $[M_\ell: M_{\ell-1}]$  (see \cite[Lem. 3.3]{SgrPG}).
We denote by 
$\boxed{\phi_{\ell, \R}}: M'_{\ell-1, \R} \to M_{\ell, \R}$
the epimorphism of real vector spaces
defined by \eqref{eq:phi-ell}. 
The dual homomorphism $\boxed{\phi_\ell^*} : N_\ell \to N'_{\ell-1}$ is
injective. It is the restriction to $N_\ell$ of the monomorphism 
\begin{equation} \label{eq:phi-ell-3}
\boxed{\phi_{\ell, \R}^*}: N_{\ell, \R} \to N'_{\ell-1, \R}, \quad
u = \sum_{i=1}^d u_i \epsilon_i \mapsto  \sum_{i=1}^d u_i \epsilon_i+ 
\langle u , 
\a_\ell - \a_{\ell-1}  \rangle  \epsilon_{d+1}  = u +  
\langle u , \a_\ell - \a_{\ell-1} \rangle  \epsilon_{d+1} .
\end{equation}

\medskip

By definition we have that $N_\ell = \{ \nu \in N_{\ell -1}  \mid \langle \nu, \a_\ell \rangle \in \Z \}$.
It follows that: 
 \begin{equation} \label{eq:rhoj}
 \phi_\ell ^* (\rho \cap N_\ell) = \bar{\r}_\ell \cap N_{\ell-1}' ,
 \end{equation}
 (see  \cite[Lemma 17]{GPRQo}). 
 Equation \eqref{eq:rhoj}, means that the pair $(\rho, N_\ell)$ and $ (\bar{\r}_\ell,  (N_{\ell-1}')_{\bar{\r}_\ell})$ are isomorphic, where  we recall that 
  $(N_{\ell-1}')_{\bar{\r}_\ell}$ is the sublattice of $N_{\ell-1}'$ spanned by $\bar{\r}_\ell \cap N_{\ell-1}'$.
 Then, we deduce an isomorphism of semigroups: 
 \begin{equation} \label{eq:rhoj-2}
 \bar{\r}_\ell^\vee
 \cap M'_{\ell-1} \longrightarrow    (\r^\vee \cap M_{\ell}) \times \ker (\phi_\ell),
 \end{equation}
 which depends on a choice of a supplementary sublattice of  $\ker (\phi_\ell) $ in $M_{\ell-1}'$
 (see \cite[Section 2.1]{FTV}).  
Notice that 
 $\ker (\phi_\ell)  $ is isomorphic to $\Z$.
We deduce from \eqref{eq:rhoj-2} the 
isomorphism  \eqref{eq:rhoj-3}. \end{proof}

\section{Toroidal embedded resolutions of a quasi-ordinary hypersurface}   \label{resolution}

In this section we outline the construction of a
\textit{toroidal embedded resolution} of a q.o. hypersurface following \cite{GPRQo, MMFGG}. 

\medskip 

We keep notations of Section \ref{qoh}. We consider an irreducible germ of quasi-ordinary hypersurface
 $(\X,o) \subset (\C^{d+1}, 0)$.  
Recall that we have fixed a complete sequence of semi-roots 
$\x_{1}, \dots, \x_{d+g+1} $ (see Definitions \ref{def:semi-root} and \ref{def: semi}).  We take  $(x_1, \dots, x_{d+1})$ as  local system of coordinates at the origin 
of $\C^{d+1}$ (see Remark \ref{rem:y}). We denote by $f \in \C\{\x \} [ \x_{d+1}]$ the quasi-ordinary polynomial 
defining $(\X,o)$.
 
 \medskip 
 
 We will define below a sequence of pairs  $(Z_\ell, o_\ell)$, where $Z_\ell$ is a normal variety of dimension $d+1$, and $o_\ell \in Z_\ell$ is a closed point, 
 together with modifications 
 \[
\psi_\ell: Z_\ell \longrightarrow Z_{\ell-1}, 
 \]
 for $\ell = 0, \dots, g$, where $Z_{0} = \C^{d+1}$ . 
By using the sequence $ \x_1, \dots, \dots, \x_{d+g+1} $ we will show that the local ring 
 $\mathcal{O}_{Z_\ell, o_\ell}$  of germs of analytic functions is isomorphic to 
 the  toric ring of convergent power series  
 \[
 \boxed{R_\ell} := \C \{ \varrho^\vee \cap M_\ell' \},
 \] 
 where by definition 
 \[ \boxed{R_0}  : =  \C \{ \x, x_{d+1} \} =  \C \{\varrho^\vee \cap M_0'\}.\]
  For $\ell \in \{0, \dots,g \}$ there exists an isomorphism
 \begin{equation} \label{eq:isoRele}
 \mathcal{O}_{Z_\ell, o_\ell} \longrightarrow  R_\ell, 
 \end{equation}
 which will provide a local toric structure of $Z_\ell$ at the point $o_\ell$.  
 We denote by $\boxed{\Psi_{\ell}} $ the composition
\[
  \Psi_{\ell} := \psi_{1} \circ \ldots \circ \psi_{\ell} : Z_\ell \to Z_0.
\]
   and by convenience we set
 $ \boxed{\Psi_0}$ the identity map of $\C^{d+1}$. 
 If $Y \subset \C^{d+1}$ is a reduced hypersurface we denote by $\boxed{Y^{(\ell)}}$ 
the \textit{strict transform} of $Y$ by $\Psi_\ell$, which is the closure of 
the preimage by $\Psi_\ell$ of the complement in $Y$ of the discriminant locus of $\Psi_\ell$.

\begin{notation} \label{def:sf} $\,$  
\begin{itemize}
\item If $h \in R_\ell$  for $\ell \in \{0,  1, \dots, g\}$,
we denote by $
\boxed{\mathcal{N}_\ell (h)}
$
the Newton polyhedron of $h$. We denote 
by $\boxed{\Sigma_\ell (h)}$ the dual subdivision of $\varrho$ defined by 
the $\mathcal{N}_\ell (h)$ and by  $\boxed{\Phi_{h}^{(\ell)}}: \varrho \to \R$ the support function 
of 
$\mathcal{N}_\ell (h)$.  We denote by $\boxed{\psi_\ell (h)}$ the Newton modification of $Z_\ell$ 
defined by the subdivision $\Sigma_\ell (h)$ (see notation \ref{newton}). 

\item If $h \in R_0$ and $\ell \in \{0,  1, \dots, g\}$,
we abuse slightly of notation by identifying $h \circ \Psi_\ell$ with its power series expansion in the ring 
$R_{\ell}$ by using the isomorphim \eqref{eq:isoRele}. 
We  denote by 
$
\boxed{\mathcal{N}_\ell (h)}
$
the Newton polyhedron of $h \circ \Psi_\ell$ in the ring $R_{\ell}$. 
We denote by $\boxed{\Sigma_\ell (h)}$ the subdivision of $\varrho$ defined by 
the $\mathcal{N}_\ell (h)$ and by  $\boxed{\Phi_{h}^{(\ell)}}: \varrho \to \R$ the support function 
of 
$\mathcal{N}_\ell (h)$.  
We denote by $\boxed{\psi_\ell (h)}$ the Newton modification $\psi (h \circ  \Psi_\ell)$ of $Z_\ell$ 
defined by the subdivision $\Sigma_\ell (h)$. 

\item
We denote by $\boxed{K}$ the quotient field of the ring ${R}_0$. 
If $h =  \frac{a}{b} \in K$ with $a, b \in R_0 \setminus \{ 0 \}$, 
 we set $\boxed{\Phi_{h}^{(\ell)} }:= \Phi_{a}^{(\ell)} - \Phi_{b}^{(\ell)}$. 
 \end{itemize}
\end{notation}

Let us explain first how to define the map $\psi_1$. 
By definition, $\boxed{y_0}  :=\x_{d+1}$ is a semi-root of degree one of $\x_{d+g+1} = f$ 
(see Remark \ref{rem:y}).  Then, 
$y_0$ is also a semi-root of degree one $\x_{d+j+1}$ 
for $j \in \{1 , \dots, g\}$. 
This implies that the power series expansion of $\x_{d+j+1}$ in $R_0$  is of the form 
\begin{equation} \label{forma-f}
\x_{d+j +1} = (y_0^{n_1} - c_1 \, \x^{n_1 \alpha_1})^{n_2 \cdots n_j} + \cdots
\end{equation}
 where $c_1 \in \C^*$ and  the exponents of the terms which are not written 
do not belong to the compact edge  of the Newton polyhedron $\Newton_0 (x_{d+j+1}) \subset
\varrho$, for $j\in \{1, \dots, g\}$.
The modification 
\[
\boxed{\psi_{1}} : \boxed{Z_{1}} \longrightarrow Z_{0},
\]
 is the Newton modification $\psi (x_{d+g+1})$ (see Definition \ref{def-Newtonm}).
By \eqref{forma-f}, the Newton fan $\Sigma (x_{d+g+1})$  coincides with the fan 
$\Sigma_1 $, which was introduced in
 terms of the characteristic exponents in Definition \ref{def:SigmaL}.
By \eqref{forma-f}, we have that $\Newton_0 (x_{d+j+1}) = n_2 \cdots n_j  \Newton_0 (x_{d+2})$. This implies that the Newton fan $\Sigma (x_{d+j+1})$ is also equal to $\Sigma_1$,  for $j\in \{1, \dots, g\}$.

\begin{lemma} {\cite[Lemma 18]{GPRQo}}. \label{lem:intersection}
The intersection of the strict transform 
 $\X_{d+j+1}^{(1)}$
with the exceptional fiber $\psi_1^{-1} (0)$  
is a point  $o_{1}$  
counted with multiplicity 
$n_2 \cdots n_j $,  for  $j \in \{1, \dots, g\}$.
\end{lemma}
\begin{proof}
We have explained above that the Newton fan $\Sigma (x_{d+j+1})$ is equal to $\Sigma_1$.
By Lemma \ref{lem-0-dim},
 the strict transform $\X_{d+j+1}^{(1)}$ does not intersect the $0$-dimensional orbits of $Z_1$. 
Then,  the intersection of $\X_{d+j+1}^{(1)}$ with the exceptional fiber $\psi_1^{-1} (0)$ is contained in 
the orbit   $\orb_{\bar{\rho}_{1}}$ of the $d$-dimensional cone $\bar{\rho}_{1} \in \Sigma_1$.
The orbit $\orb_{\bar{\rho}_{1}}$ is a one-dimensional torus embedded as a closed subset in the chart 
$Z_{\bar{\rho}_{1} , N_{0}'} \subset Z_{1}$. 
The function 
\begin{align*}
w_{1} \coloneqq y_0^{n_{1}} x^{- n_{1} \a_{1}} \in \C[\bar{\rho}_1^{\perp} \cap M_0' ]  \subset \C [\bar{\rho}^\vee \cap M_0'], 
\end{align*}
is a unit on the chart $Z_{\bar{\rho}_{1} , N_{0}'} \subset Z_{1}$, and the ring  $\C [w_{1}^{\pm 1}] = \C[\bar{\rho}_1^{\perp} \cap M_0' ] $ is  isomorphic to 
the coordinate ring of the orbit $\orb_{\bar{\rho}_{1}}$.
It follows from \eqref{forma-f} that in the chart $Z_{ \bar{\rho}_{1}, N_0'}$, 
we can factor:
\begin{align} \label{eq:form1-factor}
x_{d+j+1}  \circ \psi_{1}
= \x^{n_1 \cdots n_j  \a_{1}} 
\cdot \boxed{x_{d+j+1}^{(1)}}, 
\end{align}
and then $\X_{d+j+1}^{(1)} = Z( x_{d+j+1}^{(1)})$.
By \eqref{forma-f}, the restriction of $x_{d+j+1}^{(1)}$ to the orbit  $\orb_{\bar{\rho}_{1}}$ is 
\begin{equation} \label{eq: form1j}
(x_{d+j+1}^{(1)})_{|  \orb_{\bar{\rho}_{1}} } =  (w_{1} - c_{1})^{n_{2} \cdots n_j}.
\end{equation}
This implies that  the intersection of  $\X_{d+j+1}^{(1)}$ with $\psi_1^{-1} (0)$ is 
equal to the point $\boxed{o_1}$  defined $w_{1} -  c_{1}  = 0$,
counted with multiplicity one if $j=1$ or $n_2 \cdots n_j$ if $j>1$. 
\end{proof}

\medskip
We have the following description of the local toric structure of $Z_1$ at the point $o_1$: 
\begin{lemma} \label{lem:ring2}
The ring of germs of holomorphic functions of the variety $Z_1 $ at the point $o_1$ is isomorphic to 
$R_1 = \C \{ \varrho^\vee \cap M_1' \}$. 
\end{lemma}
\begin{proof}
By \eqref{eq:rhoj-3} the coordinate ring of the affine chart  $Z_{{\bar{\r}},N_{0}'}$
is isomorphic to 
$\C[ \rho^\vee \cap M_1] [w_1^{\pm 1}]$, and the orbit $\orb_{\bar{\r}_1}$ is defined by 
the vanishing of $x^{m}$ for $m \in  \rho^\vee \cap M_1$.
By \eqref{eq: form1j} we have that 
\[
\boxed{y_1}: = x_{d+2}^{(1)}  = w_1 -c_1 + \dots 
\]
where the terms which are not written vanish on the orbit of $\bar{\r}_1$. 
This implies that the germ of $Z( x_{d+2}^{(1)} )$ at $o_1$ is isomorphic to the germ  $Z_{\r, N_1}$ at the origin.
The maximal ideal of 
the local ring $\mathcal{O}_{Z_1, o_1}$ is generated by 
 $y_1$ together with monomials
  $\x^{\b_1}, \dots, \x^{\b_m}$,
 where $\b_1, \dots, \b_m$ generate the semigroup $\r^\vee \cap M_1$.
Taking into account \eqref{eq:rep-local} and notation \ref{not:basic}, we have an isomorphism: 
\begin{equation} \label{eq:iso-local1}
\mathcal{O}_{Z_1, o_1}  = \C\{{\r}^\vee  \cap M_1\} 
\{  y_1 \} \longrightarrow  \C\{\r^\vee  \cap M_1\} \{  x^{\varepsilon_{d+1}} \}  =   \C\{ \varrho^\vee \cap M_1' \}  = R_1 
\end{equation}
which fix $\C\{{\r}^\vee  \cap M_1\}$ and which sends $y_1$ to $x^{\varepsilon_{d+1}}$.
\end{proof}
We identify the local rings $R_1$ with $\mathcal{O}_{Z_1, o_1}$ by using 
the isomorphism \eqref{eq:iso-local1}. 

\begin{proposition} {\cite[Proposition 19]{GPRQo}} \label{tqo-key2}  We have the following properties for $ j \in \{2, \dots, g\}$.

\begin{enumerate}
\item 
The series $ \x_{d+j+1}^{(1)} \in R_1$ is irreducible and quasi-ordinary with respect the base $Z_{\r, N_1}$, 
has characteristic exponents $
		\a_{2} - \a_{1},   \ldots,  \a_{j} - \a_{1}$,
			and characteristic integers 
			$
			n_{2} , \ldots , \ n_{j}$.
\item \label{item 2 tqo} $ \x_{d+2}^{(1)},  \dots, \x_{d+j+1}^{(1)}$ is a  sequence of characteristic semi-roots of 
$ \x_{d+j+1}^{(1)}$.
\end{enumerate}
\end{proposition}

\medskip 

This process can be iterated for $\ell \in \{ 1, \dots, g-1 \}$ starting at $\ell = 1$,
which was the previous case. 
For $\ell \in \{ 2 , \dots, g-1  \}$ and  $j \in \{ \ell +1, \dots, g\}$ we have defined a function  
$\boxed{y_{\ell}} := x_{d + \ell +1}^{(\ell)}$, which is a defining function of the germ of strict transform ${\X_{d +\ell +1}^{(\ell)}}\subset Z_{\ell}$  
at the point $\boxed{o_{\ell}} \in Z_\ell$. 
We have an isomorphism of local rings 
\begin{equation}
\mathcal{O}_{Z_\ell, o_\ell} = \C \{ \rho^\vee \cap M_{\ell}\} \{ y_{\ell} \}   \longrightarrow  \C \{ \varrho^\vee \cap M_{\ell} ' \}  =  \C \{ \r^\vee \cap M_{\ell} \} \{{ x^{\varepsilon_{d+1}}} \} = R_\ell, 
\end{equation}
which fixes $\C \{ \rho^\vee \cap M_{\ell}\}$ and sends $y_{\ell}$ to $x^{\varepsilon_{d+1}}$. 
We have defined functions $ \x_{d+\ell +1}^{(\ell)},  \dots, \x_{d+j+1}^{(\ell)}$ in $R_\ell$ defining the strict transforms 
of $\X_{d+\ell +1}^{(\ell)}, \dots , \X_{d+j+1}^{(\ell)} \subset Z_\ell$, respectively. 
These functions verify the following property:
\begin{lemma}\cite[Lemma 5.8]{MMFGG} \label{lem:Manuel}
We can factor in the local ring $R_\ell$:
\[
x_{d+j+1} \circ \Psi_\ell = \x^{n_{\ell} \cdots n_j \g_\ell} \cdot x_{d+j+1}^{(\ell)}.
\] 
In particular, if $j = \ell$ we obtain that 
\begin{equation} \label{eq: yelle}
x_{d+\ell+1} \circ \Psi_\ell = \x^{n_{\ell}  \g_\ell} \cdot x_{d+\ell+1}^{(\ell)} =  \x^{n_{\ell}  \g_\ell} \cdot y_\ell.
\end{equation}
\end{lemma}

\medskip

The key fact is that Proposition \ref{tqo-key2}, can be generalized to the case when we start from 
a quasi-ordinary power with respect to a toric base.
By  {\cite[Proposition 19]{GPRQo}}, the series $ \x_{d+j+1}^{(\ell)} \in R_\ell$ is irreducible and quasi-ordinary with respect  to the base $Z_{\r, N_\ell}$, 
has characteristic exponents $
		\a_{\ell+1} - \a_{\ell},   \ldots,  \a_{j} - \a_{1}$,
			and characteristic integers 
			$
			n_{\ell+1} , \ldots , \ n_{j}$.
			In addition,  $ \x_{d+\ell +1}^{(\ell)},  \dots, \x_{d+j+1}^{(\ell)}$ is a sequence of characteristic semi-roots of 
$ \x_{d+j+1}^{(\ell)}$. 
{This implies that we have a power series expansion of $ x_{d+j+1}^{(\ell)}  \in R_\ell$ of the form:
\begin{equation} \label{forma-f-ell}
x_{d+j+1}^{(\ell)} =  
 \left( y_{\ell}^{n_{\ell+1}} - c_{\ell+1}\, \x^{n_{\ell+1} ( \alpha_{\ell +1} - \a_\ell) }	\right)^{n_{\ell +2} \cdots n_j} + \cdots,
\end{equation}
where $c_{\ell +1} \in \C^*$, 
and  the terms which are not written 
in \eqref{forma-f-ell} do not lie on the compact edge of the  Newton polyhedron
$\Newton_\ell (x_{d+j+1} )\subset
\varrho^\vee$.

\begin{remark} \label{rem: twovertices}
By \eqref{forma-f} and \eqref{forma-f-ell}, the Newton polyhedron $ \Newton_\ell (x_{d+j+1} )$ has only one compact edge with vertices 
\begin{equation} \label{eq: twovertices}
n_{\ell+1} \cdots n_j \varepsilon_{d+1} \mbox{ and } n_{\ell+1} \cdots n_j (\a_{\ell+1} - \a_\ell),
\end{equation}
 hence 
\begin{equation} \label{eq:homothetic}
 \Newton_\ell (x_{d+j+1} ) = n_{\ell+2} \cdots n_j  \, \Newton_\ell (x_{d+\ell+2} ), 
 \end{equation}
for $\ell \in \{ 0 , \dots, g-1  \}$ and  $j \in \{ \ell +1, \dots, g\}$. 
\end{remark}

We denote by  
\[\boxed{\psi_{\ell +1}} : \boxed{Z_{\ell +1}} \longrightarrow Z_{\ell},
\]
 the Newton modification $\psi_\ell (x_{d+g+1})$ (see Definition \ref{def-Newtonm} and notation \ref{def:sf}).
This modification is defined in terms of the Newton fan $ \Sigma_\ell (x_{d+g+1})$, which it is equal to the fan $\Sigma_\ell$ introduced in terms of the characteristic exponents in Definition \ref{def:SigmaL}. 

\medskip 
On the chart $Z_{\bar{\r}_{\ell+1}, N_\ell'} \subset Z_{\ell +1}$ 
we can factor 
\begin{equation} \label{form:fact-ell}
x_{d+j+1}^{(\ell)}   \circ \psi_{\ell +1}
= \x^{- n_{\ell+1} \cdots n_j ( \alpha_{\ell +1} - \a_\ell) }
\cdot \boxed{ x_{d+j+1}^{(\ell +1)}}. 
\end{equation}
The strict transform $ \X_{d+j+1}^{(\ell +1)} $ is equal to $Z( x_{d+j+1}^{(\ell +1)}) $ and intersects 
the exceptional fiber $\psi_{\ell +1}^{-1} (o_\ell)$ at a point $o_{\ell +1}$ with multiplicity 
$n_1 \cdots n_j / n_1\cdots n_\ell$. 

\medskip

We set $\boxed{y_{\ell +1}} := x_{d +\ell +2}^{(\ell +1)}$. 
We obtain then an isomorphism of local rings 
\begin{equation}
\mathcal{O}_{Z_{\ell+1}, o_{\ell+1}} = \C \{ \rho^\vee \cap M_{\ell+1}\} \{ y_{\ell+1} \}   \longrightarrow  \C \{ \varrho^\vee \cap M_{\ell+1} ' \}  =  \C \{ \r^\vee \cap M_{\ell+1} \} \{x^{\varepsilon_{d+1}}\} = R_{\ell+1}, 
\end{equation}
which fixes $\C \{ \rho^\vee \cap M_{\ell+1}\}$ and sends $y_{\ell+1}$ to $x^{\varepsilon_{d+1}}$. 

\medskip 
At the end of this process we obtain that $\boxed{y_{g}} = x_{d + g+1}^{(g)} \in R_g$ is a defining function of the strict transform 
$\X_{d+g+1}^{(g)}$ of the quasi-ordinary hypersurface $\X_{d+g+1} =\X$.  We have the following result:

\begin{theorem} \cite[Th. 1]{GPRQo}
The proper morphism $\Psi_g = \psi_{1} \circ \cdots \circ \psi_{g}$
is an embedded normalization of the quasi-ordinary hypersurface germ $D \subset \C^{d+1}$. 
\end{theorem}

Notice that the Newton polyhedron $\Newton_g (x_{d+g+1} )\subset \varrho$ has only one vertex, which is equal to $\varepsilon_{d+1}$.
The associated dual fan $\Sigma_{g} (x_{d+g+1})$
 is the fan $\Sigma_{d+g+1}$ of the lattice $N_{g}'$  consisting of the faces of cone $\varrho$  (see Definition \ref{def:SigmaL}).

\medskip 
We have the following description of  the toroidal embedding associated with the partial resolution $\Psi_g$ of the quasi-ordinary hypersurface $\X$ (see Section \ref{sec-TE}).

\begin{proposition}  \label{prop:iso-complex} \cite{GPRQo}
The pair  consisting of the
normal variety $Z_g$ and the complement $\boxed{U} \subset Z_g$ of the
total transform of the divisor $D_1 +\cdots + D_{d+g+1}$  by the modification 
$\Psi_{g}$, defines a toroidal embedding without
self-intersection. 
The associated  \textit{conic polyhedral complex $\boxed{\bar{\Theta}}$
 with integral structure}  
is obtained from the fans $\Sigma_j$, for $j=1, \dots, g+1$, by 
\[{\bar{\Theta}} = (  \bigsqcup_{j=1}^{g+1} \Sigma_j )/_\sim,
\]
where the relation $\sim$  identifies the lattice cone
$\r $ of the fan $\Sigma_\ell$  with respect to the lattice $N_\ell$, 
 with the cone  $ \phi_{\ell, \R}^* ( \r) = \overline{\r}_{\ell -1}$ of the fan $\Sigma_{\ell -1}$ with respect to the  lattice $N_{\ell -1}$. 
\end{proposition}
Recall that the homomorphism $\phi_{\ell, \R}^*$ in Proposition \ref{prop:iso-complex}
was defined in \eqref{eq:phi-ell-3}, for 
 $\ell \in \{ 2, \dots, g +1 \}$.}

\medskip 

We consider a special kind of regular subdivisions $\Sigma_{\ell}^{\textrm{reg}}$ of the fans $\Sigma_{\ell}$ for $\ell \in \{ 1, \dots, g \}$. 
\begin{definition} \label{def: compatible}
 The  regular subdivisions $\Sigma_{\ell}^{\textrm{reg}}$ of the fans $\Sigma_{\ell}$ are \textit{compatible} if
for every cone $\theta \in \Sigma_{\ell -1}^{\textrm{reg}}$ such that $\theta \subset \r$, the cone $\phi_{\ell, \R}^* ( \theta) \subset \r_\ell$ belongs to 
$\Sigma_{\ell}^{\textrm{reg}}$, for $\ell \in \{ 2, \dots, g+1 \}$.
\end{definition}
We build compatible subdivisions by 
starting with  an arbitrary regular subdivision $\Sigma_{1}^{\textrm{reg}}$, and then we define a regular subdivision 
$\Sigma_{2}^{\textrm{reg}}$ taking into account the compatibility condition imposed by $\Sigma_{1}^{\textrm{reg}}$, and 
then repeating the argument inductively.
If the regular subdivisions  $\Sigma_{\ell}^{\textrm{reg}}$,  $\ell \in \{ 1, \dots, g \}$, are compatible,  then
we obtain a regular subdivision $\boxed{\bar{\Theta}^{\textrm{reg}}}$ of the complex 
$\bar{\Theta}$ by glueing the lattice cones $\theta \in \Sigma_{\ell}^{\textrm{reg}}$ with $\theta \subset \r$, with  $\phi_{\ell, \R}^* ( \theta) \in \Sigma_{\ell-1}^{\textrm{reg}}$, for $\ell \in \{1, \dots, g+1\}$. 
Then, we have a modification 
\begin{equation} \label{eq:tor-res}
 \boxed{\Psi_{\mathrm{reg}}} : \boxed{Z} \longrightarrow Z_g
 \end{equation} 
  associated with a regular subdivision $\bar{\Theta}_{\mathrm{reg}}$ of the
conic polyhedral complex $\bar{\Theta}$, which is a resolution of singularities of $Z_g$ (see Section \ref{sec-TE}). We have the following theorem:

\begin{theorem}\label{P4thmres} \cite[Th. 1]{GPRQo}
The composition $\boxed{\Psi} :=  \Psi_{\mathrm{reg}} \circ \Psi_g : Z  \longrightarrow Z_0$ is 
a \textit{toroidal embedded resolution} of the quasi-ordinary  hypersurface germ $(\X,o)$. 
\end{theorem}

\begin{remark} \label{rem: emb-res}
The map $\Psi$ is an embedded resolution of $D$ according to the definition given at the introduction of this paper. Notice that 
$\Psi$ is not necessarily an isomorphism over $Z_0 \setminus \mathrm{Sing} (\X)$ (see \cite{McEwan-Nemethi}).
This is why we apply a particular definition of embedded resolution. 
We have that 
$\Psi$ induces an isomorphism $Z \setminus \Psi^{-1}(\X) \rightarrow Z_0 \setminus \X$.
This is consequence of \cite{McEwan-Nemethi}, see more generally \cite[Lem. 22]{GPRQo}. 
If $\X'$ is the strict transform of $\X$ we have that $\Psi$ induces an isomorphism over 
$\X \setminus \mathrm{Sing} (\X)$. This can be seen by using that 
$\Psi_g$ is an embedded normalization and that the regular subdivisions considered 
do not subdivide the regular cones.
\end{remark}

\begin{definition} \label{rem: iso-complex2}
   The boundary divisor of the toroidal embedding associated with  $\Psi_g$ is equal to 
 \begin{equation} \label{eq:bdZg}
 \boxed{\partial Z_g}=  \sum_{\theta \in \bar{\Theta}(1)} D_\theta,
 \end{equation}
  where $\bar{\Theta}(1)$ denotes the one dimensional skeleton of the the complex $\bar{\Theta}$,  and $D_\theta$ denotes
 the divisor associated to the ray $\theta$ (see Section \ref{sec-TE}).
 The boundary divisor of the toroidal embedding associated with  $\Psi$ is 
  \begin{equation} \label{eq:bdZ}
  \boxed{\partial Z}=  \sum_{\theta \in \bar{\Theta}^{\mathrm{reg}}(1)} D_\theta.
  \end{equation}  
\end{definition}

 \section{Quasi-monomial valuations of the toroidal embedded resolution} \label{sec-qmv}

In this section we describe a set of valuations of $\C^{d+1}$ associated with the toroidal embedded resolutions of the quasi-ordinary hypersurface. Let us start with the definition of semivaluation.

\begin{definition}  
   \label{def:semival}
   Denote by $\boxed{\bar{\R}_{\geq 0}}$ the set $\R_{\geq 0 } \cup \{ \infty \}$. 
    Let $A$ be a commutative ring. A \textit{semivaluation of 
    $A$} is a function $\nu : A \to {\bar{\R}_{\geq 0}} $ such that   $\nu(1) =0$ and $\nu(0) = + \infty$ and 
    for every $h, h' \in A$ we have:
          \begin{equation} \label{eq: def-semi-val}
          \begin{array}{lcl}
          \nu(h h') & = &  \nu(h) + \nu(h'),
          \\
           \nu(h+h') & \geq & \min \{ \nu(h) , \nu(h') \}.
           \end{array}
           \end{equation}
\end{definition}
A semivaluation $\nu$ is a  \textit{valuation} if $0$ is the only element $h$ of $A$ such that $\nu(h) = + \infty$.

\medskip

If $ h \in R_0$, and if $\ell \in \{1, \dots,g\}$,  we consider 
$h \circ \Psi_\ell$ as an element of  the ring $R_{\ell}$ in the following definition (see notation \ref{def:sf}). 
\begin{definition} \label{def: U} 
If $u \in \varrho$ and if $	\ell \in \{ 0, \dots, g\}$ we denote by $\boxed{\nu_{u}^{(\ell)}}$ the valuation 
defined by:
\begin{equation}
\nu_{u}^{( \ell)} ( h )  := \nu_u (h \circ \Psi_\ell),   \mbox{ for } h \in R_0. 
\end{equation}
Recall that $\nu_u$ is the monomial valuation associated to the vector $u$ (see \eqref{eq:monomial valuation}).
The valuations of the form $\nu_u^{(\ell)}$ are the
  \textit{quasi-monomial} valuations associated with the toroidal resolution process.
 \end{definition}
By notations  \ref{def:sf}
and  \ref{newton} we have that:
\begin{equation} \label{eq: Phi-h}
 \nu_{u }^{(\ell)} (h) = \Phi_h^{(\ell)} (u). 
\end{equation}

In the case of plane curve singularities, quasi-monomial valuations are studied by many authors including 
Spivakovsky \cite{S} and Favre and Jonsson \cite{TVT}.

\begin{remark} \label{rem:exc-div} Take  $\ell \in \{0, \dots, g\} $, $u \in \varrho$ and consider
the quasi-monomial valuation $\nu_{u}^{(\ell)}$. 
We can write $u =v  + r {\epsilon}_{d+1}$ where 
$v = \sum_{i = 1}^d \nu_{u}^{(\ell)}  (x_i) {\epsilon}_i \in \rho$ and 
$r = \nu_u^{(\ell)}  ( y_{\ell} )$. This means that the valuation $\nu_{u}^{(\ell)}$ is determined by its values on 
the monomials $x_1, \dots, x_d, y_{\ell}$ in the ring 
$R_{\ell}$.
\end{remark}
\begin{remark} \label{def:exc-div}
If $\ell \in \{0, \dots, g\}$ and $u \in \varrho  \cap N_\ell'$ is a primitive integral vector for the lattice $N_\ell'$,  the valuation $\nu_{u}^{(\ell)}$ is the \textit{vanishing order valuation} associated with the toroidal divisor 
$\boxed{D_{u}^{(\ell)}}$.
In particular, if $\ell \in \{0, \dots, g\}$ and $i \in \{1, \dots, d+1\}$ we denote by $\boxed{\epsilon^{(\ell)}_i}$ the primitive integral vector 
 of the lattice $N_{\ell}'$ which belongs to the face $\R_{\geq 0}  \epsilon_i$ of the cone $\varrho$. 
 We denote the divisor ${D_{\epsilon_{i}^{(\ell)}}^{(\ell)}}$ associated with the valuation
 $\nu_{\epsilon_{i}^{(\ell)}}^{(\ell)}$
simply by 
$\boxed{D_{i}^{(\ell)}}$.

\medskip 

If $u$ does not belong to a ray of the cone $\varrho$, then the valuation $\nu_{u}^{(\ell)}$ 
is the \textit{divisorial valuation} associated with the toroidal exceptional divisor 
$D_{u}^{(\ell)}$.   This divisor  appears in a model $Z_{\ell}'$ obtained from $Z_{\ell}$ by a toroidal modification defined by a subdivision of $\Sigma_\ell$ containing the ray spanned by $u$ 
 (see Section \ref{sec-toric}).

\medskip 
If $u$ belongs to a ray of the cone $\varrho$, then  the divisor $D_{u}^{(\ell)}$ appears on $Z_{\ell}$ and it may be non-exceptional.
In particular, 
if $\ell \in \{ 0,1, \dots, g \}$ the valuation 
$\nu_{\epsilon_{d+1}}^{(\ell)}$ is the \textit{vanishing order valuation} 
along $\X_{d+ \ell +1}$, that is, 
if $h \in R_0$ and we factor $h = x_{d +\ell +1}^k \cdot h'$ for some integer $k \in \N$ 
and some $h' \in R_0$ which is not divisible by $h$,  then $\nu_{\epsilon_{d+1}}^{(\ell)} (h) = k$. 
This valuation is not centered at the maximal ideal of $R_0$.  
In this case,  the associated divisor $D_{\epsilon_{d+1}}^{(\ell)} = D_{d+1}^{(\ell)}$ 
is the strict transform of $\X_{d +\ell + 1}$ on $Z_\ell$. 
Similarly, we have  the vanishing order valuations
$\nu_{\epsilon_{i}}^{(0)}$ 
 along $\X_i = D_{\epsilon_{i}}^{(0)}$, for $i \in \{1,\dots, d+1\}$.
 \end{remark}

 \begin{remark} \label{rem: iso-complex}
By the construction of the toroidal normalization process, if $	\ell \in \{ 2, \dots, g\}$, $u  \in \rho \subset \varrho$ and $\tilde{u} := \phi_{\ell, \R}^* (u) \in N_{\ell -1, \R}'$ (see \eqref{eq:rhoj}), then  the quasi-monomial valuations 
$\nu^{(\ell)}_{u}$ and $\nu_{\tilde{u}}^{(\ell-1)}$ coincide. 
In particular, if $u \in \rho \cap N_{\ell}$ is a primitive vector, then 
the exceptional divisors $D_{u}^{(\ell)}$ and $D_{\tilde{u}}^{(\ell -1)}$ coincide. 
\end{remark}

 \begin{remark} \label{rem: iso-complex3}
With the notation introduced in remark \ref{def:exc-div} and  Definition \ref{rem: iso-complex2}, we have that  the components of the boundary divisor 
 $\partial Z_g$ 
 are $D_i^{(\ell)}$ for $\ell \in \{ 0, \dots, g\}$ and $i\in \{1, \dots,d\}$, where we may have that 
 $D_i^{(\ell)} = D_i^{(\ell')}$ for some $i \in \{1, \dots, d \}$ and $\ell, \ell ' \in \{0, \dots, g\}$ (see Example 
 \ref{Ex:Two pairs qo jumping numbers}).
 \end{remark}

 \subsection{Generating sequences of the quasi-monomial valuations}

The notion of generating sequence of a valuation is classical in valuation theory  (see for instance \cite{S} 
in the case of two dimensional regular local rings).
In this section we prove that $x_1, \dots, x_{d+g+1}$ is a generating sequence of 
any quasi-monomial valuation associated to the toroidal resolution process.

\begin{definition}
Let $\nu$ be a semivaluation of $R_0$. 
Let $z_1, \dots, z_m \in R_0$ be a system of generators of the maximal ideal of $R_0$. 
We say that  $z_1, \dots, z_m $ is a \textit{generating sequence} of 
$\nu$ if for any $a \in \R_{\geq 0}$ the set 
 $\boxed{\mathcal{P}_\nu(a)}  : = \{ h \in R_0 \mid \nu (h) > a \} $
 is an ideal of $R_0$ generated by monomials in $z_1, \dots, z_m$. 
\end{definition}

 We introduce now some useful lemmas about  support functions of the transforms of the semi-roots in 
 the embedded normalization process (see notation \ref{def:sf}).
 
 \begin{lemma} \label{lem: t0}
If $\ell \in \{ 0, \dots, g \}$,  then we have that 
$\Phi_{x_i}^{(\ell)} = {\varepsilon}_i$, for  $i \in \{1, \dots, d \}$, and 
$\Phi^{(\ell)}_{y_{\ell}} = {\varepsilon}_{d+1}$. 
\end{lemma}
 \begin{proof}
 This follows by applying Definition   \ref{def:sf}  taking into account Remark \ref{newton}.
\end{proof}

\begin{lemma} \label{lem: t2}
If  $\ell \in \{1, \dots,g \}$ and $j \in \{ 0, \dots, \ell -1 \}$,  then we have that 
\begin{equation}
\Phi^{(\ell)}_{x_{d+j+1}} = \g_{j+1}.
\end{equation}
\end{lemma}
\begin{proof}
Take $j \in \{ 0, \dots, \ell -1\}$. 
Then, the function 
$w_{j+1} = y_j^{n_{j+1}} \cdot \x^{-n_{j+1} (\a_{j+1} - \a_j)}$ is a unit  in the ring
$\C \{ \varrho^\vee \cap M_{j+1} '\}$. It is also a unit when we consider it 
in the ring $\C \{ \varrho^\vee \cap M_{\ell} '\}$. 
Then:
\begin{equation}
\Phi^{(\ell)}_{y_{j}} 
= \Phi^{(\ell)}_{ \x^{ \a_{j+1} - \a_j} }  =  \a_{j+1} - \a_j.
\end{equation}
Since 
$x_{d+j+1}  \circ \Psi_j = \x^{n_j \g_j } \cdot y_j$ on $\C \{ \varrho^\vee \cap M_{j} ' \}$,  we obtain that 
\[
\Phi^{(\ell)}_{x_{d+j+1}} =  \Phi^{(\ell)}_{\x^{n_j \g_j }} +   \Phi^{(\ell)}_{y_{j}} =  n_j \g_j + \a_{j+1} - \a_j \stackrel{\eqref{rel-semi}}{=} \g_{j+1}. 
\]
\end{proof}

\begin{remark} By notation \ref{not:basic}, we consider  $\r^\vee$ as a face of $\varrho^\vee$. Hence, we can see  the semigroup 
$\Gamma$ as a subset of  $\varrho^\vee$.  
 By Lemma \ref{lem: t0} and remark \ref{rem:exc-div}, if 
$\g \in \Gamma_\ell $,  then $\x^\gamma \in R_\ell$ and 
$\Phi_{\x^\gamma}^{(\ell)} = \g$. 
\end{remark}

\begin{lemma} \label{lem: t1} $\,$

\begin{enumerate}[label=(\alph*)]
\item \label{t1a}
If  $\ell \in \{0, 1, \dots,g \}$ and $j \in \{ \ell, \dots,g \}$,  then we have that 
\begin{equation} \label{eq: t1a}
\Phi_{x_{d+j+1}}^{(\ell)} = n_{\ell} \cdots n_j \g_\ell  + \Phi_{x_{d+j+1}^{(\ell)}}^{(\ell)}
\end{equation}
In addition, 
\begin{equation}\label{eq: t1c}
\Phi_{x_{d+\ell+1}}^{(\ell)} = n_{\ell}  \g_\ell  + \varepsilon_{d+1}. 
\end{equation}

\item  \label{t1b}  If $\ell \in \{0, \dots, g-1 \}$  and if $j \in \{ \ell +1, \dots,  \leq g \}$,  then:
\begin{equation} \label{eq: t1b}
\Phi_{x_{d+\ell+2}^{(\ell)}}^{(\ell)}    =  \min\{     n_{\ell+1}  (\a_{\ell +1} - \a_\ell), \,   n_{\ell+1}  \varepsilon_{d+1} \},
\quad \mbox{ and } \quad 
\Phi_{x_{d+j+ 1 }^{(\ell)}}^{(\ell)} = n_{\ell+2} \cdots n_{j} \, \Phi_{x_{d+\ell+2 } ^{(\ell)}}^{(\ell)}.
\end{equation} 
\end{enumerate}
\end{lemma}
\begin{proof}
In the case $\ell = 0$, the assertion follows from Formula \eqref{forma-f} taking into account that by definition $\g_0 = 0$ and  $n_0 =1$.  
Assume then that $\ell > 0$. 
The first statement follows by Lemma  \ref{lem:Manuel}, and by applying  Definition \ref{def:sf}.
The formula \eqref{eq: t1b} is  is a reformulation in terms of support functions of  \eqref{forma-f-ell}. 
\end{proof}

\begin{proposition} \label{prop: genseq1} Let us take $\ell \in  \{0, \dots, g\}$ and $u \in \varrho$. 
Then, $x_1, \dots, x_{d + \ell +1}$ is a generating sequence of the valuation $\nu_{u}^{(\ell)}$. 
\end{proposition}
\begin{proof}
Take $h \in R_0, h \ne 0$. Its $(x_1, \dots, x_{d+\ell+1})$-expansion is of the form
\begin{equation} \label{eq: hexp}
h = {\sum}_{I} 
c_{I} \, 
x_1^{i_1}  \cdots   \x_{d+\ell+1}^{i_{d+\ell+1}}, 
\end{equation}
where if $c_I \ne 0$, $I = (i_1, \dots, i_{d+ \ell +1} )$, then
\begin{equation} \label{eq: cota-I}
0\leq i_{d+j}< n_{j}, \ \mathrm{for}\ j \in \{ 1,\ldots ,\ell \}. 
\end{equation}
(see Definition \ref{def: expansion}).
Let us denote 
$
\mathcal{M}_I := x_1^{i_1}  \cdots   \x_{d+\ell+1}^{i_{d+\ell+1}}.
$
By Lemmas \ref{lem: t0}, \ref{lem: t2} and  formula \eqref{eq: t1c} 
the support function  $ \Phi_{\mathcal{M}_I}^{(\ell)} $ is equal to 
\[
\Phi_{\mathcal{M}_I}^{(\ell)}= \sum_{j=1}^d i_j \varepsilon_j + \sum_{j =0}^{\ell -1} i_{d+j+1}  \g_{j+1} + i_{d + \ell + i}  
( n_\ell \g_\ell + \varepsilon_{d+1}) \in M_{\ell}'.
\]
Let $\mathcal{M}_I$ and $\mathcal{M}_{I'}$ appear in the expansion \eqref{eq: hexp}.
If $\Phi_{\mathcal{M}_I}^{(\ell)}
= \Phi_{\mathcal{M}_{I'}}^{(\ell)}$ then representing this vector in terms of the basis $\{ \varepsilon_i \}_{i=1}^{d+1}$
we get that $ i_{d + \ell + i}   =  i_{d + \ell + i} '$. 
Then, 
\begin{equation} \label{eq: n(I)}
 \Phi_{\mathcal{M}_I}^{(\ell)}  -  i_{d + \ell + i} (n_\ell \g_\ell + \varepsilon_{d+1})  =  \Phi_{\mathcal{M}_{I'}}^{(\ell)}  -  
 i_{d + \ell + i} ' (n_\ell \g_\ell + \varepsilon_{d+1}) \in \Gamma_\ell.
\end{equation}
This equality implies that  $I = I'$ by  \eqref{eq: cota-I}  and 
by the distinguished representation of the elements of the semigroup $\Gamma_\ell$ given in 
Lemma \ref{lem:unique-sg} . 

By \eqref{eq: Phi-h}
we have that 
$\nu_{u}^{(\ell)} (\mathcal{M}_I) = \Phi_{\mathcal{M}_I}^{(\ell)} (u) $ and $\nu_{u}^{(\ell)} (h) = \Phi_{h}^{(\ell)} (u) $. 
By the previous discussion we get that 
\[  \nu_{u}^{(\ell)} (h) = \min_{c_I \ne 0} \nu_{u}^{(\ell)} (\mathcal{M}_I).\] 
Therefore, for any real number $a$ the valuation  ideal,
 $\mathcal{P}_{\nu_{u}^{(\ell)}} (a) $ is generated by monomials in 
 $x_1, \dots, x_{d+\ell+1}$, that is,  these functions define a
 generating sequence of $\nu_{u}^{(\ell)}$. 
\end{proof}

As any finite subset containing a generating sequence is a generating sequence we obtain:
\begin{corollary}
For any $\ell \in  \{0, \dots, g\}$ and $u \in \varrho$ we have that 
$x_1, \dots, x_{d + g +1}$ is a generating sequence of the valuation $\nu_{u}^{(\ell)}$. 
\end{corollary}

Next, we prove that under certain hypothesis on  
$u$ and $\ell$,  we can 
compute the value $\nu_{u}^{(\ell)} (h)$ from the values of 
the monomials in the
$(x_1, \dots, x_{d+g+1})$-expansion of $h$, for  any $h \in R$. 

\begin{proposition} \label{prop: genseq2} 
Take $\ell \in \{0, \dots, g-1\}$ and  $u \in \bar{\sigma}_{\ell +1}^-$ (see Definition \ref{def:SigmaL}).
Consider the  $(x_1, \dots, x_{d+g+1})$-expansion 
\eqref{Form:Expansion of any germ in terms of semi-roots} of a nonzero function 
$h \in R_0 $. 
Then, we have that 
\[
\nu_{u}^{(\ell)} (h) = \min_{c_I \ne 0}  \nu_{u}^{(\ell)} ( x_1^{i_1} \dots x_{d+g+1}^{i_{d+g+1} }).
\]
\end{proposition}
\begin{proof}
Take a term $\mathcal{M}_I$ appearing in the expansion \eqref{Form:Expansion of any germ in terms of semi-roots}
with nonzero coefficient $c_I$. 
By Lemmas \ref{lem: t0}, \ref{lem: t2} and \ref{lem: t1} the support function 
$\Phi_{\mathcal{M}_I}^{(\ell)} $ is of the form: 
\[
\begin{array}{lcl}
\Phi_{\mathcal{M}_I}^{(\ell)} & = &
 \sum_{j=1}^d i_j \varepsilon_j + \sum_{j =0}^{\ell -1} i_{d+j+1}  \g_{j+1} + i_{d + \ell + 1}  
( n_\ell \g_\ell + \varepsilon_{d+1})
\\
&  & + 
(i_{d+\ell+2} n_{\ell+1} + \cdots + i_{d+g+1} n_{\ell+1} \cdots n_g)
 \min\{      \a_{\ell +1} - \a_\ell, \,   \varepsilon_{d+1} \}.
 \end{array}
\]
Let us set
\begin{equation} \label{eq: mr(I)}
\left\{
\begin{array}{lcl}
r(I)  & :=  &  i_{d + \ell + 1}   + i_{d+\ell+2} n_{\ell+1} + \cdots + i_{d+g+1} n_{\ell+1} \cdots n_g  \in \N
\\
m(I) & := &  \sum_{j=1}^d i_j \varepsilon_j + \sum_{j =0}^{\ell -1} i_{d+j+1}  \g_{j+1} + i_{d + \ell + 1}  
n_\ell \g_\ell +  r(I) \varepsilon_{d+1} \in M_{\ell}'.
\end{array}
\right.
\end{equation}
By Lemmas \ref{lem: t0}, \ref{lem: t2} and \eqref{eq: t1c}, the support function $\Phi_{x_i}^{(\ell)}$ is linear, for $i\in \{1, \dots, d+\ell +1\}$. This means  that the Newton polyhedron $\mathcal{N}_\ell(x_i)$ 
has only one vertex, for $i\in \{1, \dots, d+\ell +1\}$. By Remark \ref{rem: twovertices}, 
the Newton polyhedron $ \Newton_\ell (x_{d+j+1} )$
has two vertices \eqref{eq: twovertices}, and it is homothetic to  $\Newton_\ell (x_{d+g+1} )$, for $j \in \{ \ell +1, \dots, g\}$.
 It follows that
the Newton polyhedron
$\mathcal{N}_\ell (\mathcal{M}_I)$ has at most two vertices, and one of them is $m(I)$. 
If $\mathcal{F}_I$ is the compact face of maximal dimension of $\mathcal{N}_\ell (\mathcal{M}_I)$ then 
we get  $\mathcal{N}_\ell (\mathcal{M}_I)  = \mathcal{F}_I + \varrho^\vee$.
If $\dim (\mathcal{F}_I) =1 $ then  $\mathcal{F}_I$ is parallel to the vector 
$ \varepsilon_{d+1} -  \a_{\ell +1} +\a_\ell$. 
In particular, we have that  if $ (m_1', \dots, m'_{d+1})$ are the coordinates of a vector $m' \in \mathcal{F}_I$ with respect to the basis $\{ \varepsilon_i \}_{i=1}^{d+1}$, then 
\begin{equation} \label{eq: aux-fu}
m'_{d+1} \leq r(I) \mbox{ with equality only when } m' = m(I).
\end{equation}

By \eqref{eq: Phi-h} we have that 
\begin{equation} \label{eq: sfell}
\nu_{u}^{(\ell)} ( {\mathcal{M}_I} ) = \Phi_{\mathcal{M}_I}^{(\ell)} (u) =
\langle u, m(I) \rangle,
\end{equation}
where the last equality follows from the hypothesis $u \in \bar{\sigma}_{\ell +1}^-$. 

\medskip

\noindent
\textit{Claim}: If $c_I, c_{I'} \ne 0$ and  if $m(I) = m({I'})$ then $I = I'$. 

\noindent
Let us prove this claim first. 
Since 
$\g_1, \dots, \g_\ell$ belong to the subspace spanned by $\{ \varepsilon_i \}_{i=1}^{d}$, 
and $\{ \varepsilon_i \}_{i=1}^{d+1}$ is a basis of  $M_{\ell, \R}'$, 
the hypothesis $m(I) = m({I'})$ imply that $r(I) = r({I'})$. 
By using the elementary argument appearing in the proof of \cite[Lemma 7.2]{PPPDuke} 
we obtain that 
\begin{equation} \label{eq: equal-I}
i_{d + j + 1} = i_{d+j +1}', \mbox{ for } j \in \{\ell, \dots, g\}.
\end{equation}
The  hypothesis $m(I) = m({I'})$ and the equality \eqref{eq: equal-I}  implies that
\[
\sum_{j=1}^d i_j \varepsilon_j + \sum_{j =0}^{\ell -1} i_{d+j+1}  \g_{j+1} = \sum_{j=1}^d i_j ' \varepsilon_j + \sum_{j =0}^{\ell -1} i_{d+j+1} ' \g_{j+1} \in \Gamma_\ell.
\]
Then, the claim follows 
by  condition \eqref{eq: cota-I}  and 
by the distinguished representation of the elements of the semigroup $\Gamma_\ell$  of
Lemma \ref{lem:unique-sg}.  
This ends the proof of the claim.

\medskip 
Let us set 
\[
\mathcal{S} := \{ I \mid c_I \ne 0,  \nu_{u}^{(\ell)} ( {\mathcal{M}_I} )  
= 
\min_{c_{I'} \ne 0}\{  \nu_{u}^{(\ell)} ( \mathcal{M}_{I'} ) \} \},  
\]
Take $I_0 \in \mathcal{S}$ such that $r(I_0) = \max_{I \in \mathcal{S}}  r(I)$.
Since $\nu_u^{(\ell)}$ is a valuation,  and ${\mathcal{M}_{I_0}}$ is a term of $h$, 
 we have that $ \nu_{u}^{(\ell)} ( {\mathcal{M}_{I_0}} )   \leq  \nu_{u}^{(\ell)} ( h)$.  
It  is  enough to prove that $m(I_0) \in
\mathcal{N}_\ell (h)$,  since this implies that $ \nu_{u}^{(\ell)} ( h)   \, = \min \{ \langle u, v \rangle \mid v \in \mathcal{N}_\ell (h) \} \leq  \langle 
u, m(I_0) \rangle  \stackrel{\eqref{eq: sfell}}{=}  \nu_{u}^{(\ell)} ( {\mathcal{M}_{I_0}} )$.

\medskip 
If $m(I_0) \in \mathcal{N}_\ell  (\mathcal{M}_{I_1})$,  then we get that 
 \begin{equation} \label{eq:inequ}
  \nu_{u}^{(\ell)} ( {\mathcal{M}_{I_1}} )  \stackrel{\eqref{eq: sfell}}{=} \langle u, m(I_1) \rangle \leq 
 \nu_{u}^{(\ell)} ( {\mathcal{M}_{I_0}} ) 
  \stackrel{\eqref{eq: sfell}}{=}  \langle u, m(I_0) \rangle.
\end{equation}
In addition, equality holds since ${I_0} \in \mathcal{S}$, hence  $I_1 \in \mathcal{S}$. 
The equality in \eqref{eq:inequ} also implies that $m(I_0)$ and $m(I_1)$ belong  to  the compact face 
$\mathcal{G}_u$ of $\mathcal{N}_\ell ({\mathcal{M}_{I_1}})$ defined by $u$.  Then, we distinguish two cases:

- If $\mathcal{G}_u$ is a vertex then we have $m(I_1)= m(I_0)$. 

- Otherwise, $\mathcal{G}_u$ is equal to the compact edge $\mathcal{F}_{I_1}$ of $\mathcal{N}_\ell ({\mathcal{M}_{I_1}})$.  By \eqref{eq: aux-fu} applied to $m(I_0)$, we get that $r(I_0) \leq r(I_1)$. Since $I_1\in \mathcal{S}$, the hypothesis $r(I_0) = \max_{I \in \mathcal{S}}  r(I)$, implies that $r(I_1) \leq r(I_0)$, thus $r(I_1) = r(I_0)$. 
By the case of equality in \eqref{eq: aux-fu}, applied to $m(I_0)$, it follows that $m(I_1)= m(I_0)$. 

In both cases,  $m(I_1)= m(I_0)$ holds. By the claim we get $I_0 = I_1$. 
This implies that $m(I_0)$ is a vertex of  $ \mathcal{N}_\ell (h)$. 
\end{proof}

\section{Log-discrepancies of exceptional divisors} 
\label{Sect:log-discrepancy vector}

 In this section we fix a toroidal embedded resolution \eqref{eq:tor-res} of the quasi-ordinary hypersurface $\X$. 
 We keep the notations of the previous sections. 
We describe below the log-discrepancies of the  divisors
$D_{u}^{(\ell)}$ associated to a primitive vector $u \in \varrho \cap M_\ell'$ 
for some $\ell \in \{ 0, \dots, g\}$ (see remark \ref{def:exc-div}).

\medskip 
Recall from Section \ref{Sec: jn} the definition of the relative canonical divisor 
of a modification between two smooth complex varieties
and the log-discrepancies of its exceptional prime components.

\begin{definition}
We denote 
\begin{equation} \label{eq: lambda_l}
\boxed{\lambda_0} :=  \varepsilon_1 + \cdots + \varepsilon_{d+1}  \in M_0' \mbox{ and } 
\boxed{\lambda_\ell} := \alpha_\ell + \lambda_0  \in M_\ell', \mbox{ for } \ell \in \{ 1, \dots, g\} . 
\end{equation}
We say that $\lambda_\ell$ is the \textit{log-discrepancy vector of depth} $\ell \in \{ 0, 1, \dots, g \}$ in
the toroidal resolution.  
\end{definition}

\begin{lemma} \label{lem: equal-l}
If $v \in \rho$ and $\ell \in \{ 1, \dots, g\}$,  then
$
\langle v, \l_\ell \rangle = \langle \phi_\ell^* (v), \l_{\ell -1} \rangle, 
$
where $\phi_\ell^*$ is described in \eqref{eq:phi-ell-3}.
\end{lemma}
\begin{proof}
We apply the definitions using that   $\langle \epsilon_{d+1} , \a_{\ell -1} \rangle =0$ and  
$\langle \epsilon_{d+1} , \l_0 \rangle =1$: 
\[
\begin{array}{lcl}
\langle \phi_\ell^* (v), \l_{\ell -1} \rangle &  \stackrel{\eqref{eq:phi-ell-3}}{=} &  \langle v  + \langle v, \a_\ell - \a_{\ell-1} \rangle \epsilon_{d+1} , \lambda_{\ell-1} \rangle 
=  \langle   v  + \langle v, \a_\ell - \a_{\ell-1} \rangle \epsilon_{d+1} , \a_{\ell -1} + \l_0 \rangle
\\
& = & \langle v, \a_{\ell -1} \rangle + \langle v, \l_0 \rangle 
+ \langle v, \a_\ell - \a_{\ell-1} \rangle  (  \langle \epsilon_{d+1} , \a_{\ell -1} \rangle +  \langle \epsilon_{d+1} , \l_0 \rangle )
\\
& = & 
  \langle v, \a_{\ell -1} \rangle + \langle v, \l_0 \rangle 
+ \langle v, \a_\ell - \a_{\ell-1} \rangle
= \langle v, \a_\ell  + \l_0 \rangle
= \langle v, \l_\ell \rangle. 
\end{array}
\] 
\end{proof}

The following result expresses the log-discrepancies of the divisors  $D_{u}^{(\ell)}$
appearing in the toroidal embedded resolution in terms of the characteristic exponents of the 
quasi-ordinary hypersurface (compare it with
Proposition 8.2 and Remark 8.4 of
\cite{MMFGG}). 
Proposition \ref{prop:log-discrepancy vector}
is a generalization of a result of
Favre and Jonsson, \cite[Proposition D.1]{TVT}, 
see also \cite[Theorem 8.18 and Proposition 8.26]{GGP19} and 
\cite[Section7]{JonssonDynamicsBerkovich}, which state a similar result in the case of plane curve singularities.

\begin{proposition} \label{prop:log-discrepancy vector}
Take $\ell \in \{ 0, \dots, g \}$.  Let  $u \in \varrho \cap N_\ell'$ be a primitive vector.
Then, the log-discrepancy $\lambda_{D_{u}^{(\ell)}}$ of the prime divisor $D_{u}^{(\ell)}$ is equal to 
$ \langle u , \l_\ell \rangle$. 
\end{proposition}
\begin{proof}

We prove the statement by induction on $\ell$.
If $\ell = 0$, we take
 vectors $u , a_2,  \dots, a_{d+1} \in \varrho \cap N_0'$ 
defining a basis of the lattice $N_0'$. 
We obtain a monomial map $\psi: \C^{d+1} \to \C^{d+1}$ given by:
\[
x_{1}=z_{1}^{u_{1}}z_{2}^{a_{2,1}}\ldots z_{d+1}^{a_{d+1,1}},  \quad \dots,  \quad  x_{d+1}=z_{1}^{u_{d+1}}z_{2}^{a_{2,d+1}}\ldots z_{d+1}^{a_{d+1,d+1}} .
\]
The divisor $D_{u}^{(0)}$ is defined by $z_1 = 0$.
Then, we get that 
\[
\psi^{*} (dx_{1} \wedge \ldots \wedge dx_{d+1}) = z_{1}^{\langle u,  \l_0 \rangle -1} \ldots 
	z_{d+1}^{\left\langle a_{d+1}, \l_0  \right\rangle -1} 
	dz_{1} \wedge \ldots \wedge dz_{d+1}.
\]
Therefore,  $\l_{D_{u}^{(0)}} = \langle u , \l_0 \rangle$. 
 Assume that the statement is true for $\ell -1$. 
 
 \medskip 
We take a basis $u_1, \dots, u_{d+1}$ of the lattice $N_\ell'$ such that $u_{d+1} = \epsilon_{d+1}$, $u_1, \dots, u_d \in \r$ and $u \in \R_{\geq 0} u_1 + \cdots + \R_{\geq 0} u_{d+1}$ 
(see notation \ref{not:basic}). Then, the cone 
$\theta := \R_{\geq 0} u_1 + \cdots + \R_{\geq 0} u_d$ is regular for the lattice $N_\ell$. We may assume that 
$\theta$ belongs to the regular subdivision $\Sigma_{\ell}^{\mathrm{reg}}$ and 
$\bar{\theta} := \phi_\ell^* (\theta) \subset \bar{\r}_{\ell-1} \in \Sigma_{\ell-1}^{\mathrm{reg}}$ 
by the compatibility condition in  Definition \ref{def: compatible}. 
Consider the lattice homomorphism $\phi_{\ell}^*: N_\ell \to N_{\ell -1} '$ defined by \eqref{eq:phi-ell-3}.
By construction of the toroidal resolution the  divisors,
$D_{u_i}^{(\ell)}$ and $D_{\phi^*(u_i)}^{(\ell-1)}$ are the same for $i \in \{1, \dots, d \}$
 (see Remark \ref{rem: iso-complex}).

On the chart $Z_{\bar{\theta}} \subset Z_{\ell}^{\mathrm{reg}}$ 
corresponding to  $\bar{\theta}$,  we take local coordinates 
$(s_1, \dots, s_{d+1})$ at the point 
where the orbit of $\bar{\theta}$ intersects the strict transform  $D_{\epsilon_{d+1}}^{(\ell)}$ of $\X_{d+\ell+1}$. 
We may assume that 
\[
Z(s_i) = D_{\phi^*(u_i)}^{(\ell-1)}, \mbox{ for } i=1, \dots, d,  \mbox{ and } Z(s_{d+1}) = 
D_{\epsilon_{d+1}}^{(\ell)}.
\]

By the induction hypothesis we have that: 
\begin{equation}\label{eq: indhyp}
\lambda_{D_{\phi^*(u_i)}^{(\ell-1)}} = \langle \phi^*(u_i) ,  \l_{\ell -1},\rangle.
\end{equation}
By Lemma \ref{lem: equal-l} we obtain:
 \begin{equation} \label{eq: lambda-b} 
b_i := \lambda_{D_{u_i}^{(\ell)}} = \langle u_i ,  \l_{\ell} \rangle, \mbox{ for } i \in \{ 1, \dots, d \}. 
\end{equation}
Since the divisor $D_{\epsilon_{d+1}}^{(\ell)}$ is the strict transform of $\X_{d+\ell +1}$,
its  log-discrepancy is equal to one, and  
\begin{equation} \label{eq: lambda-b2}
 b_{d+1} : = \l_{D_{\epsilon_{d+1}}^{(\ell)}} = \langle u_{d+1}, \l_{\ell}  \rangle  = 1. 
\end{equation}

Next, we take a basis $w_1, \dots, w_{d+1}$  of $N_{\ell}'$ such that 
\begin{equation} \label{eq:incl-condition}
w_1 = u, \quad   \mbox{ and } 
 \s := \R_{\geq 0} w_1 + \cdots + \R_{\geq 0} w_{d+1} \subset \theta + \R_{\geq} \epsilon_{d+1}.
 \end{equation}
 We may assume that the cone $\sigma$ belongs to a regular subdivision $\Sigma_{\ell}^\mathrm{reg}$.
By \eqref{eq:incl-condition} we expand
\begin{equation} \label{eq:base}
w_i = \sum_{j=1}^{d+1} w_{i,j} u_j, \mbox{ for } i =1,\dots, d+1
\end{equation} where  $w_{i,j} \geq 0$. 
We consider the monomial map $\psi_\sigma: \C^{d+1} \to \C^{d+1}$ defined by 
\[
s_{1}=z_{1}^{w_{1,1}} z_{2}^{w_{2,1}}\ldots z_{d+1}^{w_{d+1,1}},  \quad \dots ,  \quad 
s_{d+1}=z_{1}^{w_{1, d+1}}z_{2}^{w_{2,d+1}}\ldots z_{d+1}^{w_{d+1,d+1}} .
\]
The divisor $D_u^{(\ell)}$ appears on this chart as $Z(z_1)$, thus by \eqref{eq: lambda-b}  and \eqref{eq: lambda-b2} its
log-discrepancy $\l_{D_u^{(\ell)}}$ is one plus the order of vanishing of 
\[
\psi_\s^* (s_1^{b_1-1} \cdots s_{d+1}^{b_{d+1}-1} ds_1 \wedge \cdots  \wedge ds_{d+1} )
\] on the divisor $Z (z_1)$.
We write $b := (b_1, \dots, b_{d+1})$  and $w_i := (w_{i,1}, \dots, w_{i, d+1})$ for $i=1, \dots, d+1$, and we compute: 
\[
\psi^*_\s ( s_1^{b_1-1} \cdots s_{d+1}^{b_{d+1}-1} ) =  \psi^*_\s ( s_1^{b_1} \cdots s_{d+1}^{b_{d+1}} ).
(\psi^*_\s ( s_1 \cdots s_{d+1} ))^{-1} = 
z_1^{\langle w_1, b \rangle} \cdots  z_{d+1}^{\langle w_{d+1} , b \rangle} (\psi^*_\s ( s_1 \cdots s_{d+1} ))^{-1}.
\]
We can reformulate the computation done in the case $\ell =0$ to get:
\[
\psi_\s^* (ds_1 \wedge \cdots  \wedge ds_{d+1} ) = (\psi^*_\s ( s_1 \cdots s_{d+1} )) (z_1 \dots z_{d+1})^{-1}  
dz_1 \wedge \cdots  \wedge dz_{d+1}.
\]
Therefore:
\[
\psi_\s^* (s_1^{b_1-1} \cdots s_{d+1}^{b_{d+1}-1} ds_1 \wedge \cdots  \wedge ds_{d+1} ) = 
z_1^{\langle w_1, b \rangle} \cdots  z_{d+1}^{\langle w_{d+1} , b \rangle} (z_1 \dots z_{d+1})^{-1} dz_1 \wedge \cdots  \wedge dz_{d+1}.
\]
It follows that $\l_{D_u^{(\ell)}} = {\langle w_1, b \rangle}$. 
By \eqref{eq: lambda-b}, \eqref{eq: lambda-b2}  and \eqref{eq:base} we deduce that
\[
\l_{D_u^{(\ell)}} = {\langle w_1, b \rangle} = \sum_{j=1}^{d+1} w_{1, j} b_j  = 
\sum_{j=1}^{d+1}  w_{1, j} \langle u_i,  \l_{\ell}  \rangle 
=    \langle  \sum_{j=1}^{d+1} w_{1, j} u_i, \l_{\ell} \rangle 
= \langle u, \l_\ell \rangle. 
\]
\end{proof}

\begin{remark}The notion of log-discrepancy has been extended to valuation spaces of normal varieties 
by 
Boucksom, de Fernex, Favre and Urbinati in \cite{BdFF15}. 
\end{remark}

\section{Multiplier ideals of a quasi-ordinary hypersurface} \label{sec: miqh}

In this section we give a description of the multiplier ideals of an irreducible germ of quasi-ordinary hypersurface. 
We consider a toroidal log-resolution $\Psi: Z \to \C^{d+1}$ of the quasi-ordinary hypersurface germ $(\X, 0) \subset (\C^{d+1}, 0)$
described in Section  \ref{resolution} in terms of a  regular subdivision $\bar{\Theta}^{\textrm{reg}}$ of the conic polyhedral 
complex $\bar{\Theta}$. The boundary divisor ${\partial Z}$ of the associated 
toroidal embedding  has components $D_\theta$ for $\theta$ running through the one skeleton
$\bar{\Theta}^{\mathrm{reg}}(1)$ of the 
complex $\bar{\Theta}^{\textrm{reg}}$ (see  \eqref{eq:bdZ}). The component  $D_\theta \subset Z$ is either
 exceptional or the strict transform of one of the hypersurfaces 
$\X_i$, for $i \in \{1, \dots, d+g+1\}$.
Recall that  $x_{d+g+1}$ is a defining function of the q.o. hypersurface $\X = \X_{d+g+1}$. We have that 
\[ \Psi^* (x_{d+g+1} ) = \mathcal{O}_Z (- \sum_{\theta \in \bar{\Theta}^{\mathrm{reg}}(1)}\nu_{D_\theta} (x_{d+g+1} )  D_\theta ). \] 

\medskip 

The following lemma is an application of 
 the valuative description  \eqref{eq:form MI} of the multiplier ideals.
\begin{lemma} \label{lem:mi-theta}
Take $\xi \in \Q_{\geq 0}$. 
A function $h \in  \mathcal{O}_{\C^{d+1},0} $ belongs to the multiplier ideal $\mathcal{J}( \xi \X)_0$ 
if and only if  the valuative inequality 
\begin{equation} \label{eq-mi-gen}
 \nu_{D_\theta} (h) > \xi  \, \nu_{D_\theta} (x_{d+g+1}) - \lambda_{D_\theta}. 
 \end{equation}
holds for $D_\theta$ in the support of the boundary divisor $\partial Z$.
In addition, if $0 < \xi < 1$,  then $h \in \mathcal{J}( \xi \X)_0$ if and only if
condition \eqref{eq-mi-gen} holds for every exceptional divisor
$D_\theta$ in the support of the boundary divisor $\partial Z$.
\end{lemma}
\begin{proof}
Notice that if $D_\theta$ is the strict transform of $\X_{d+g+1}$, then $\nu_{D_\theta} (x_{d+g+1} ) = 1$ while if 
 $D_\theta$ is the strict transform of $\X_i$, for some $i \in \{1, \dots, d+g\}$,  then $\nu_{D_\theta} (x_{d+g+1} ) = 0$.  
If $D_\theta$ is not exceptional  its log-discrepancy $\l_{D_\theta}$ is equal to one.
The inequality  \eqref{eq-mi-gen} holds since either 
 $\nu_{D_\theta} (x_{d+g+1}) = 1$ if $D_\theta = D_{d+g+1}$ 
   or $\nu_{D_\theta} (x_{d+g+1}) = 0$ otherwise.
In both cases  the assertion follows since  $0 <\xi < 1$.
\end{proof}

We  prove now that some of the conditions \eqref{eq-mi-gen} are superfluous.
The following theorem is a generalization of Theorem 4.9 and Corollary 4.17 of \cite{GGGR24}, which was inspired by  Smith and Thompson  \cite{STIrrelExcDiv} in the case of plane curve singularities.

\begin{theorem} \label{Thm: reduction} 
Take $\xi \in \Q_{\geq 0}$. 
The following conditions are equivalent: 
\begin{enumerate}[label=(\alph*)]
\item  \label{a}
The function $h \in  \mathcal{O}_{\C^{d+1},0}$ belongs to the multiplier ideal 
$\mathcal{J}(\xi \X)_0$.

\item \label{b}
The polyhedron $ \Newton_{\ell -1} (h) + \lambda_{\ell -1}$ is contained in the interior of  the polyhedron  $ \xi \Newton_{\ell -1}(x_{d+g+1})$, for 
$\ell \in \{1, \dots, g+1\}$.

\item   \label{c}
The valuative inequality \eqref{eq-mi-gen}
holds for 
 $D_\theta$ in the support of the boundary divisor $\partial Z_g$ of the toroidal embedded normalization (see \eqref{eq:bdZg}).
\end{enumerate}

\end{theorem}
\begin{proof}
By Lemma \ref{lem:mi-theta} and the description 
of the complex $\bar{\Theta}^{\textrm{reg}}$ in terms of the regular fans ${\Sigma_{\ell}}^{\textrm{reg}} $, we get that
 $h \in \mathcal{J}(\xi \X)_0$ if and only if condition 
\eqref{eq-mi-gen} holds for every ray $\theta$ of   ${\Sigma_{\ell}}^{\textrm{reg}}$, and for $\ell \in \{1, \dots, g+1 \}$. 
Take $\ell \in \{ 1, \dots, g+1\}$. 
Let $\theta$ be a ray of the fan $\Sigma_{\ell}^{\mathrm{reg}}$. There is a primitive integral vector 
$u_\theta \in \varrho \cap N_{\ell -1}'$ such that $D_\theta = D_{u_\theta}^{(\ell -1)}$ 
(see Remark \ref{rem:exc-div}).  

\medskip
By Proposition \ref{prop:log-discrepancy vector} the log-discrepancy of $D_\theta$ is equal to:
\begin{equation} \label{eq: lambdaD} 
\lambda_{D_\theta} = \langle u_\theta, \lambda_{\ell -1} \rangle.
\end{equation}
By \eqref{eq: Phi-h}  and by definition one gets: 
\begin{equation} \label{eq: lambdaD1}
\nu_{D_\theta} (h) = \nu_{u_\theta}^{(\ell -1)} (h) = \Phi_{h}^{(\ell -1)}  (u_\theta) = \Phi_{\Newton_{\ell -1} (h)} (u_\theta),
\end{equation}
in particular
\begin{equation} \label{eq: lambdaD2}
\nu_{D_\theta} (x_{d+g+1}) =  \nu_{u_\theta}^{(\ell -1)} (x_{d+g+1})  = \Phi_{x_{d+g+1}}^{(\ell -1)} (u_\theta ) = \Phi_{ \Newton_{\ell -1}(x_{d+g+1})} (u_\theta ) .
\end{equation}
By \eqref{eq: New-xi} we obtain that 
\begin{equation} \label{eq: lambdaD4}
\xi \nu_{D_\theta} (x_{d+g+1})   = \xi \Phi_{ \Newton_{\ell -1}(x_{d+g+1})} (u_\theta)  = 
\Phi_{ \xi \Newton_{\ell -1}(x_{d+g+1})} (u_\theta) .
\end{equation}
By \eqref{eq-sf-tras} we deduce that 
\[
\Phi_{\Newton_{\ell -1} (h)}  + \lambda_{\ell -1}  =  \Phi_{\Newton_{\ell -1} (h) + \lambda_{\ell -1}}.
\]

It follows that the valuative condition \eqref{eq-mi-gen}  is equivalent to the inequality 
\begin{equation} \label{eq: lambdaD3}
 \Phi_{\Newton_{\ell -1} (h) + \lambda_{\ell -1}} (u_\theta) -  \Phi_{ \xi \Newton_{\ell -1}(x_{d+g+1})} (u_\theta) >0.
 \end{equation} 
Then, the  equivalence of conditions \ref{a}, \ref{b} and \ref{c} follows 
by Proposition \ref{prop: int-pol}.
\end{proof}

\begin{remark} \label{rem: mayor1}
Let  $\xi ' \in \Q_{>0}$ and $\xi = \xi'  -  \lfloor \xi \rfloor ' \in [0, 1)$. By Lemma \ref{Lem:periodicidad}
we have  that $\mathcal{J}(\xi' \X)_0 = (x_{d+g+1})^{\lfloor \xi' \rfloor} \cdot  \mathcal{J}(\xi \X)_0$.
By Lemma \ref{lem:mi-theta} and Theorem \ref{Thm: reduction} we get that 
$h \in \mathcal{J}(\xi \X)_0$ if and only if condition \eqref{eq-mi-gen} holds for every \textit{exceptional}
divisor $D_\theta$ in the support  of the boundary divisor $\partial Z_g$.
\end{remark}

\begin{definition} \label{def:xiM}
If $i_1, \dots, i_{d+g+1} \in \Z_{\geq 0}$ and $\mathcal{M} = x_1^{i_1} \dots x_{d+g+1}^{i_{d+g+1}}$,  we 
define the rational number 
\begin{equation}
\boxed{\xi_{\mathcal{M}}}= \min \left\{ ({\nu_{D_\theta} ( \mathcal{M} ) + \lambda_{D_\theta}})({  \nu_{D_\theta} (x_{d+g+1}) })^{-1}  \mid D_\theta \mbox{ is a prime exceptional divisor of } \partial Z_g 
 \right\}.
\end{equation}
\end{definition}

 In the following theorem we give a
 description of the generators  of the multiplier ideals of $\X$ at $0$ and of their jumping numbers. 
 It is a higher dimensional generalization of Theorem 4.20 of \cite{GGGR24}. 
\begin{theorem} \label{rem: mon-ideal}
Take a rational number $\xi \in (0,1)$.  
\begin{enumerate}[label=(\alph*)]
 \item  \label{a-M} 
 The multiplier ideal  $\mathcal{J}(\xi \X)_0$ is generated by those generalized monomials 
 $\mathcal{M}$ in $x_1, \dots, x_{d+g+1}$ such that  
 $\xi < \xi_{\mathcal{M}} \leq  \xi + 1$. 
 \item   \label{b-M} 
 A rational number $\xi$ is a jumping number of the multiplier ideals of $\X$ at $0$ if and only if there exists
 a  generalized monomial $\mathcal{M}_0$ in $x_1, \dots, x_{d+g+1}$ such that  $\xi = \xi_{\mathcal{M}_0}$. 
\end{enumerate}
\end{theorem}
\begin{proof}
\ref{a-M} 
Let $\mathcal{I}_\xi$ be the ideal of $\mathcal{O}_{\C^{d+1}, 0}$ generated by those  generalized  monomials $\mathcal{M}$ in $x_1, \dots, x_{d+g+1}$ 
such that $\xi < \xi_{\mathcal{M}} $. 
By  Theorem \ref{Thm: reduction}  \ref{c}  it follows that 
$\mathcal{I}_\xi \subset \mathcal{J}(\xi \X)_0$. 

\medskip 
Let us check the opposite inclusion.  By Theorem \ref{Thm: reduction} and Remark \ref{rem: mayor1}, one has that 
 $h \in \mathcal{J}(\xi \X)_0$ if and only if condition 
\eqref{eq-mi-gen} holds for any exceptional divisor $D_\theta$ in the support of the boundary divisor $\partial Z_g$. 
By Remark \ref{def:exc-div} and Proposition \ref{prop:iso-complex} 
we have that $D_\theta = D_{u}^{(\ell)}$ for 
$\ell \in \{ 0, \dots, g-1 \}$ and $u \in \bar{\r}_{\ell +1}$ 
defining a one dimensional face of the cone $ \bar{\r}_{\ell +1}$.
 
 \medskip
 Let us consider the $(x_1, \dots, x_{d+g+1})$-expansion 
\eqref{Form:Expansion of any germ in terms of semi-roots} of $h$ (see Lemma \ref{Lema:Expansion semi-roots}): 
\[
h = {\sum}_{I} 
c_{I} \, 
x_1^{i_1}  \cdots   \x_{d+g+1}^{i_{d+g+1}}.
\]
Since $u \in  \bar{\r}_{\ell +1}$ and $ \bar{\r}_{\ell +1} \subset   \bar{\s}_{\ell +1}^- $, 
Proposition \ref{prop: genseq2}, implies that 
\[
\nu_{D_\theta} (h)  = \min_{c_I \ne 0} \{ \nu_{D_\theta} (x_1^{i_1}  \cdots   \x_{d+g+1}^{i_{d+g+1}})  \}. 
\]
If  $c_I \ne 0$,
the inequality 
\begin{equation}   \label{eq: red-h-D2}
\nu_{D_\theta} (x_1^{i_1}  \cdots   \x_{d+g+1}^{i_{d+g+1}})  \geq \nu_{D_\theta} (h) >  \xi \nu_{D_\theta} (x_{d+g+1})  - \lambda_{D_\theta},  
\end{equation} 
holds for every exceptional prime divisor $D_\theta$ in the support of the boundary divisor $\partial Z_g$.
Then, $x_1^{i_1}  \cdots   \x_{d+g+1}^{i_{d+g+1}}  $ belongs to the multiplier ideal $\mathcal{J}(\xi \X)_0$
by Theorem \ref{Thm: reduction} and Remark \ref{rem: mayor1}.
This implies that $\mathcal{J}(\xi \X)_0 \subset \mathcal{I}_\xi$.

Assume now that  $\mathcal{M} = x_1^{i_1} \dots x_{d+g+1}^{i_{d+g+1}}$  is a generalized monomial in $\mathcal{J}(\xi \X)_0$ such that 
$\xi + 1 < \xi_{\mathcal{M}}$.  We check that $\mathcal{M}$ is not a generator of the multiplier ideal $\mathcal{J}(\xi \X)_0$. 
Since $1 < \xi_{\mathcal{M}}$, we have that 
 $\mathcal{M} $
 belongs to the ideal $\mathcal{J}(1 \cdot \X)_0 = (x_{d+g+1} ) $ by Lemma \ref{Lem:periodicidad}.
  Then, $x_{d+g+1}$ divides $\mathcal{M}$, and
 if we put $\mathcal{M}' = \mathcal{M}/ x_{d+g+1}$,  we get that 
 $\xi_{\mathcal{M}} = \xi_{\mathcal{M'}} +1$. 
By iterating  this argument we obtain a monomial $\mathcal{M}''$, dividing  $\mathcal{M}$
 and such that $\xi < \xi_{\mathcal{M}''} \leq \xi + 1$.  This ends the proof of \ref{a-M}. 
 
 \medskip

\ref{b-M}
If $\xi$ is a jumping number of the multiplier ideals of $\X$ at $0$, then there exists $0 <r$ 
small enough such that 
\[
\mathcal{J}(\xi \X)_0 \subsetneq \mathcal{J}(\xi ' \X)_0 \mbox{ and } \mathcal{J}(\xi ' \X)_0 = \mathcal{J}(\xi '' \X)_0,
\mbox{ for } \xi', \xi'' \in (\xi -r, \xi). 
\]
By the case  \ref{a-M} in Theorem \ref{rem: mon-ideal}, there exist a  generalized  monomial $\mathcal{M}_{I_0}$ in $x_1, \dots, x_{d+g+1}$ such that 
\[
\mathcal{M}_{I_0} \in  \mathcal{J}(\xi ' \X)_0 \mbox{ if } \xi' \in (\xi -r, \xi), \mbox{ and }  \mathcal{M}_{I_0} \notin  \mathcal{J}(\xi  \X)_0.
\]
By \eqref{eq: red-h-D2}, this implies that $\xi' < \xi_{\mathcal{M}_{I_0}}$,  and 
there exists an exceptional prime divisor $D_{\theta_0} $ in the support of the  boundary divisor $\partial Z_g $,  such that 
the condition
\[
\nu_{D_{\theta_0}} (\mathcal{M}_{I_0} )  >  \xi \nu_{D_{\theta_0}} (x_{d+g+1})  - \lambda_{D_{\theta_0}},
\]
is not satisfied. 
By continuity, it follows 
that $\xi  
= 
\xi_{\mathcal{M}_{I_0}}$. 

\medskip 

Finally,  let $\mathcal{M}$ be a  generalized monomial  in $x_1, \dots, x_{d+g+1}$. 
Then, $\xi_{\mathcal{M}}$ is a jumping number by Theorem 
\ref{rem: mon-ideal} \ref{a-M}, since $\mathcal{M}$ belongs to $\mathcal{J}(\xi ' \X)_0$ for $\xi' <  \xi_{\mathcal{M}} $, while $\mathcal{M} \notin \mathcal{J}(\xi_{\mathcal{M}} \X)_0$.
 \end{proof}
\section{The local tropicalization defined by the sequence  semi-roots} \label{sec: loctrop}

In this section we study the local tropicalization of $(\C^{d+1}, 0)$ defined by the complete
sequence $x_1, \dots , x_{d+g+1}$ of semi-roots of the quasi-ordinary polynomial $f$.
The theory of local tropicalizations was 
developed in detail by Popescu-Pampu and Stepanov in  \cite{PS13}, see \cite{PS25} for an 
introduction to this topic. 
 It was applied in \cite{dFGM23} to show the existence of
embedded resolutions of curves with one toric morphism, after a suitable reembedding.
This theorem was also proved in a recent paper by Cueto, Popescu-Pampu and Stepanov,  as a corollary of their results on the local tropicalizations of  some surface singularities \cite{CPS24}.

 \medskip 
 We introduce first  a  fan $\Theta$ in $\R^{d+g+1}$ associated to the characteristic exponents of 
 the q.o. polynomial $f$.

\begin{notation}  \label{cones} 
We consider the following linear maps: 
\[
\boxed{\mathcal{T}_0}: N_{0, \R} \longrightarrow \R^{d+g +1},
\quad 
\mathcal{T}_0 (\epsilon_i) =  \epsilon_i, \mbox{ for } i = 1, \dots, d,
\]
where $\{ \epsilon_i \}_{i = 1, \dots, d} $ denotes the canonical basis of $N_{0} = \Z^d$ in the source, 
and $\{ \epsilon_i \}_{i = 1, \dots, d+g+1} $ denotes the canonical basis of $\Z^{d+g+1}$ in the target, by a slight abuse of notation. 
For  $\ell \in \{1, \dots, g\}$ we set
\[
\boxed{\mathcal{T}_\ell} : N_{0, \R} \longrightarrow \R^{d+g +1}, \quad 
\mathcal{T}_\ell (\nu) = \nu + \sum_{i=1}^\ell
 \langle \nu,  \gamma_i
\rangle \epsilon_{d + i}  +  \sum_{i=\ell+1}^{g+1}    n_\ell \cdots n_{i-1} \langle \nu, \gamma_\ell \rangle \epsilon_{d + i}. 
\]
Take $\ell \in \{ 0, \dots, g \}$. The cone 
$
\boxed{\r_\ell} := \mathcal{T}_\ell (\r)$
is $d$-dimensional, simplicial  and 
rational for the lattice $\Z^{d+g +1}$. 
We introduce the following cones of dimension $d+1$:
 \[
\boxed{\sigma^+_j} :=  \rho_j + \R_{\geq 0}  \epsilon_{d+j}, \quad \boxed{\sigma^-_j} := \rho_{j-1} + \rho_j,
\quad \mbox{ for } j =1,\dots, g, \quad \mbox{ and }  \quad  \boxed{\sigma_{g+1}} := \rho_g +
\R_{\geq 0}  \epsilon_{d+g+1}.
\]
\end{notation}

\begin{definition}  \label{complex}
We set $\boxed{\Theta_j } := \{ \s_j^+, \s_j^-, \r_j \}$ for $j =1, \dots, g$ and $\Theta_{g+1} := \{ \s_{g+1} \}$.
The set $\boxed{\Theta}$ consisting of the faces of the cones in $\bigcup_{j=1}^{g+1} \Theta_j$, is a fan
associated to the q.o.~polynomial $f$.
\end{definition}

\begin{remark} The conic polyhedral complex  $\bar{\Theta}$ associated to the toroidal resolution of 
the quasi-ordinary hypersurface
 is isomorphic to the fan $\Theta$ of Definition
\ref{complex}, by an isomorphism which preserves
the integral structure (see \cite[Section 5.4]{GPRQo}).
The definition of the fan $\Theta$ was inspired by 
 the work of Lejeune-Jalabert and Reguera \cite{LR99} on sandwiched surface singularities. 
\end{remark}

In \cite{MMFGG} we described  the fan $\Theta$ in terms of the orders of contact of the tuple 
$(x_1, \dots, x_{d+g+1})$ with formal arcs. 

\medskip 

A \textit{formal arc} of $\C^{d+1}$ centered at the origin is defined by a tuple of formal power series $\varphi \in ((t)\C[[t]])^{d+1}$. 
It defines a morphism  $ \mathrm{Spec} \C [[t]] \to \C^{d+1}$.
If  the arc $\varphi$ does not factor through the hypersurface 
of equation $x_1 \dots x_{d+g+1}=0$, then
the \textit{order of contact}  of
 $\varphi$ with 
$(x_1, 
\dots, x_{d+g+1})$ is the vector:
\[
( \orden_t (x_1 \circ \varphi), \dots, \orden_t (x_{d+g+1} \circ \varphi) )  \in \Z^{d+g+1}_{>0}.
\]

\begin{theorem} \cite[Th. 7.12]{MMFGG} \label{th-arc-trop}
The support of the fan $\Theta$ is equal to the closure of the cone spanned by the orders of contact
$( \orden_t (x_1 \circ \varphi), \dots, \orden_t (x_{d+g+1} \circ \varphi)) $,  
where $\varphi$ runs through the set of  formal arcs of $\C^{d+1}$ centered at the origin and which do not
factor through the hypersurface 
of equation $x_1 \dots x_{d+g+1}=0$.
\end{theorem}

We consider the  notion of 
local tropicalization explained in \cite[\S 6]{PS25}.

\begin{definition} \label{def:finitetrop}
We denote by $\boxed{R}$ the ring $ \C [[ \x, y ]]$ where $x= (x_1, \dots, x_d)$. 
The \textit{tropicalization map} of $R$ with respect to the tuple $(x_1, \dots, x_{d+g+1})$ is 
\[
\boxed{\trop}: \mathcal{V} \to \bar{\R}^{d+g+1}_{\geq 0}, 	\quad \trop(\nu) :=   (\nu(x_1), \dots, \nu (x_{d+g +1})). 
\]
We denote by $\boxed{\mathcal{V}}$ the set of \textit{semivaluations} of the ring $R$.  
\end{definition}

Since  $x_1, \dots, x_{d+1}$ is a sequence of generators of the maximal ideal of $\mathcal{O}_{\C^{d+1}, 0}$, 
the complete sequence of semi-roots $x_1, \dots, x_{d+g+1}$ define the  following embedding:
\begin{definition} \label{def:embedding}
We consider the  embedding of 
\begin{equation} \label{eq:embedding}
(\C^{d+1},0) \hookrightarrow (Y, 0) \subset (\C^{d+g+1},0)
\end{equation} defined by the $\C$-algebra homomorphism: 
\begin{equation} \label{eq:embedding2}
\C[[ X_1, \dots, X_{d+g+1} ]] \longrightarrow \mathcal{O}_{\C^{d+1}, 0}, \quad X_i \mapsto x_i, \mbox{ for }  i \in \{1, \dots, d+g+1\}.
\end{equation}
We denote by $\mathcal{I}_Y$ the kernel of the homomorphism \eqref{eq:embedding2}.
\end{definition}
 Remark that the image of the quasi-ordinary hypersurface $\X$  by this embedding is equal to $ Y \cap \{ X_{d+g+1} = 0\}$.
 Notice that $Y$ is not contained in any of the coordinate hyperplanes of $\C^{d+g+1}$, that is, $Y$ is an \textit{interior toric subgerm} in the sense of \cite{PS25}.

\begin{definition}
The \textit{local tropicalization}  $\mathrm{Trop} (Y)$ of $Y  \subset \C^{d+g+1}$  is the closure of the  set of 
tuples of the form $(\mu(X_1), \dots, \mu(X_{d+g+1})) \in {\R}^{d+g+1}_{> 0}$,  where  $\mu$ is a semivaluation of 
$\C[[X_1, \dots, X_{d+g+1} ]]$ such that $\mu (h)  = \infty$,   for every $h \in \mathcal{I}_Y$,  and $\mu(X_i) >0$ for $i\in \{1, \dots, d+g+1 \}$. 

An \textit{arcwise weight vector}  of $Y$ is of the form $\mathrm{ord}_t (\varphi)$, where $\varphi$ is a formal 
arc of $\C^{d+g+1}$  centered at the origin, with generic point in the torus and factoring through $Y$ (see \cite[Def. 5.5]{PS25}).

 \end{definition}

\begin{corollary}
 \label{th-trop}
 The \textit{local tropicalization}  $\mathrm{Trop} (Y)$ of $Y  \subset \C^{d+g+1}$ 
  is 
equal to the support of the fan $ \Theta$. 
\end{corollary}
\begin{proof}
Let  $\mu$ be a semivaluation of 
$\C[[X_1, \dots, X_{d+g+1} ]]$ such that $\mu (h)  = \infty$, for every $h \in \mathcal{I}_Y$,  and $\mu(X_i) >0$ for $i\in \{1, \dots, d+g+1 \}$. Then,  $\nu ( h + \mathcal{I}_Y) := \mu (h)$ defines a semivaluation of the ring 
$\C[[X_1, \dots, X_{d+g+1} ]] / \mathcal{I}_Y =  R$ such that $\trop (\nu) \in \R_{>0}^{d+g+1}$, and
every semivaluation of $R$ verifying this condition can be obtained in this way. 
This implies that 
\[ 
\mathrm{Trop} (Y) = \trop (\mathcal{V}) \cap  {\R}^{d+g+1}_{>0}. 
\]

By \cite[Th. 6.2]{PS25} the local tropicalization 
$\mathrm{Trop} (Y)$ coincides with the closure of the cone spanned by the arcwise weight vectors of $Y$.
By the definition of the embedding \eqref{eq:embedding} 
the arcwise weight vectors are precisely 
the orders of contact
$( \orden_t (x_1 \circ \varphi), \dots, \orden_t (x_{d+g+1} \circ \varphi)) $,  
where $\varphi$ is a formal arc of $\C^{d+1}$ centered at the origin and which does not 
factor through the hypersurface 
of equation $x_1 \dots x_{d+g+1}=0$.
Then, the assertion follows by Theorem \ref{th-arc-trop}. 
\end{proof}

\medskip

\subsection{The multiplier ideals and the local tropicalization} \label{sec:mi-lc}

In this section we give a reformulation of our main result about the generators of the multiplier ideals in terms of the local tropicalization. 

\begin{definition} \label{xi} 
We denote by $\boxed{\epsilon_1,  \dots, \epsilon_{d+g+1}}$   the canonical basis of $\Z^{d+g+1}$ and by 
 $\boxed{\varepsilon_1,  \dots, \varepsilon_{d+g+1}} \in  \check{\Z}^{d+g+1}$ its dual basis. Let us introduce the following elements of 
$\check{\Z}^{d+g+1}$: 
\[
\boxed{\Lambda_1} : = \sum_{i=1}^{d+1} \varepsilon_i, \quad 
\mbox{ and }   \quad
\boxed{\Lambda_j}    :=    \sum_{i=1}^d \varepsilon_i + \sum_{i=1}^{j-1} (1 - n_i) \varepsilon_{d+i}  + \varepsilon_{d+j}, \  \mbox{ for } 
j =2, \dots, g+1.
\]

\end{definition}

\begin{definition}
The \emph{log-discrepancy} function 
\[ \boxed{\Lambda}: |\Theta| \to \R_{\geq 0}
\]
is defined on a vector $u$ in the support $|\Theta|$ of the fan $\Theta$ by 
\[
\Lambda (u) = \langle \Lambda_j, u \rangle \mbox{ if } u \in |\Theta_j |, 
\mbox { for some }j \in \{ 1, \dots, g+1\}.
\]
\end{definition}

\begin{lemma}
The log-discrepancy function $\Lambda$ is well-defined
and continuous. 
\end{lemma}
\begin{proof}
This is a consequence of Lemma 2.19 of \cite{MMFGG}.
\end{proof}

\begin{notation} Take $\ell \in \{ 1, \dots, g \}$ and $i \in \{ 1, \dots, d \}$.
We denote by $u_i^{(\ell)} = \sum_{j=1}^{d+g+1} u_{i,j}^{(\ell)} \epsilon_j$ the primitive integral vector of the lattice $\Z^{d+g+1}$ spanning an edge of the cone $\r_\ell$ 
such that $u_{i,i}^{(\ell)} \ne 0$,  and
 $u_{i,j}^{(\ell)} = 0$,  for $j \in \{1, \dots, d\}$  and $j \ne i$.
\end{notation} 

By Lemmas \ref{lem: t0}, \ref{lem: t2} and 
\ref{lem: t1},  the vector $u_i^{(\ell)} $ is determined by the primitive integral vector $\epsilon_i^{(\ell)}$ of the cone $\rho$ with respect to the lattice $N_\ell$. The vector $\epsilon_i^{(\ell)}$ is determined by 
 the semigroup $\Gamma$.
We illustrate this in the Example \ref{Ex:Two pairs qo jumping numbers} below. As a consequence we have that 
$\trop(\nu_{D_i^{(\ell)}}) = u_i^{(\ell)}$.
In addition, we have that the following formula for the log-discrepancy of the divisor $D_i^{(\ell)}$:
\begin{equation} \label{eq:ld-ui}
\lambda_{D_i^{(\ell)}} = \Lambda_{\ell} (u_i^{(\ell)}).
\end{equation}
Formula \eqref{eq:ld-ui} is consequence of Remark 8.4 of \cite{MMFGG}.

\begin{lemma} \label{lem:reformulation} If $I= (i_1, \dots, i_{d+g+1}) \in  \Z_{\geq 0}^{d+g+1}$ and $\mathcal{M} = x_1^{i_1} \dots x_{d+g+1}^{i_{d+g+1}}$,  
then 
\begin{equation}
\xi_{\mathcal{M}} = \min  \left\{ \frac{ \langle u_i^{(\ell)} , I \rangle + \Lambda_{\ell} ( u_i^{(\ell)})}{ i_{d+g+1} u_{i, d+g+1}^{(\ell)} } \mid i \in \{1, \dots, d\} ,\,  \ell \in \{ 1, \dots, g \},  \mbox{ and } u_i^{(\ell)} \notin \{ \epsilon_1, \dots, \epsilon_d \}  \right\}
\end{equation}
\end{lemma}
\begin{proof}
We apply the definition of $\xi_{\mathcal{M}}$ (see Definition \ref{def:xiM}), combined with \eqref{eq:ld-ui}, and the fact that 
$\trop(\nu_{D_i^{(\ell)}}) = u_i^{(\ell)}$, for $i \in \{ 1, \dots, d \}$ and $\ell \in \{ 1, \dots, g\}$.
\end{proof}

Lemma \ref{lem:reformulation} provides a reformulation of Theorem \ref{rem: mon-ideal} in terms of the primitive integral vectors  $u_i^{(\ell)}$ defining the edges of the local tropicalization $\Theta$, which are not defined by the canonical basis vectors of $\Z^{d+g+1}$.

\begin{example} \label{Ex:Two pairs qo jumping numbers}
The fractional power series 
\begin{align*}
\zeta = x_{1}^{\frac{3}{2}}x_{2}^{\frac{1}{2}} + x_{1}^{\frac{7}{4}}x_{2}^{\frac{1}{2}} . 
\end{align*} 
$\zeta$ is a quasi-ordinary branch with characteristic exponents 
$\a_{1} = \left( \frac{3}{2} , \frac{1}{2} \right)$  and $\a_{2} = \left( \frac{7}{4} , \frac{1}{2} \right)$. 
The characteristic integers are $n_{1} = 2$ and $n_{2}=4$.   
The series $\zeta$ is a root of the  quasi-ordinary polynomial 
\begin{align*}
f = (y^{2} - x_{1}^{3}x_{2})^{4} + 2 x_1^{13} x_2^4 + x_1^{14} x_2^4 - 12 x_1^{10} x_2^3 y^2 - 2 x_1^7 x_2^2 y^4 . 
\end{align*}
In  \cite[Example  11.1]{MMFGG} a similar example of irreducible quasi-ordinary surface 
with the same characteristic exponents was considered.
We have that 
\[ x_{3} = y, \quad x_4 = y^{2} - x_{1}^{3}x_{2} \quad \mbox{ and } x_5 = f  \] 
is a characteristic sequence of semi-roots. 
The semigroup of the quasi-ordinary surface is 
$\Gamma = \Z_{\geq 0} \varepsilon_1 + \Z_{\geq 0} \varepsilon_2 + \Z_{\geq 0} \g_1 + \Z_{\geq 0} \g _2 \subset \Q_{\geq 0}^2$, where 
$
\gamma_{1} = \left( \frac{3}{2} , \frac{1}{2} \right) , \mbox{ } \gamma_{2} = \left( \frac{13}{4} , 1 \right)$.

By Definition \ref{rem: iso-complex2},  the components of the boundary divisor 
 $\partial Z_g$ are $D_i^{(\ell)}$ for $\ell \in \{ 0, 1, 2\}$ and $i\in \{1, 2\}$.  In this case,
 $D_2^{(1)} = D_2^{(2)}$, since the second coordinate of $\a_2$ and $\a_1$ coincide. 
 The exceptional components of the boundary divisor $\partial Z_g$ are 
 $D_1^{(1)}, D_2^{(1)} $ and $D_1^{(2)}$.  
 By Remark \ref{def:exc-div} the divisorial valuation associated with $D_i^{(\ell)}$ is 
 the quasi-monomial valuation $\nu_{\epsilon_i^{(\ell)}}^{(\ell)}$ where $\epsilon_i^{(\ell)}$ is 
 the primitive integral vector of the lattice $N_\ell$ on the ray $\R_{\geq 0} e_i$, for $i = 1, 2$ and $\ell = 1,2$.
 By \eqref{eq:Nele} we get that 
 \[
 \epsilon_1^{(1)} = 2 \epsilon_1, \quad   \epsilon_2^{(1)} = 2 \epsilon_2, \quad  \epsilon_1^{(2)} = 4 \epsilon_1, \quad   \epsilon_2^{(2) = 2 \epsilon_2.}
 \]
By Proposition \ref{prop:log-discrepancy vector}, the log-discrepancy of $D_i^{(\ell)}$ is equal to 
$\lambda_{D_i^{(\ell)}} = \langle \epsilon_i^{(\ell)}, \l_0 + \a_\ell \rangle$, where $\l_0 = \varepsilon_1 + \varepsilon_2$. We get: 
\[
\lambda_{D_1^{(1)}} = \langle \epsilon_1^{(1)}, \l_0 + \a_1\rangle = 5, \quad 
\lambda_{D_2^{(1)}} = \langle \epsilon_2^{(1)}, \l_0 + \a_1\rangle = 3, \quad
\lambda_{D_1^{(2)}} = \langle \epsilon_1^{(2)}, \l_0 + \a_2 \rangle = 11.
\]

Next, we compute the order of vanishing of the functions $x_1, \dots, x_5$ along these divisors.
We use that the order of vanishing of $x_j$ along the divisor
$D_i^{(\ell)}$  is equal to 
\[
 \nu_{\epsilon_i^{(\ell)}}^{(\ell)} (x_j) \stackrel{\eqref{eq: Phi-h}}{=}  \Phi_{x_j}^{(\ell)} ( \epsilon_i^{(\ell)} ),
\]
for $j \in \{1, \dots, 5\}$.
Recall that $ \Phi_{x_i}^{(\ell)}$ is the support function defined in Definition \ref{def:sf}.
By Lemma \ref{lem: t0},  we get:
\[
\begin{array}{c}
 \nu_{\epsilon_1^{(1)}}^{(1)} (x_1) = \langle \epsilon_1^{(1)}, \varepsilon_1 \rangle = 2, \quad 
  \nu_{\epsilon_1^{(1)}}^{(1)} (x_2)  = \langle \epsilon_1^{(1)}, \varepsilon_2 \rangle = 0, 
  \\
   \nu_{\epsilon_2^{(1)}}^{(1)} (x_1) = \langle \epsilon_2^{(1)}, \varepsilon_1 \rangle = 0, \quad 
  \nu_{\epsilon_2^{(1)}}^{(1)} (x_2)  = \langle \epsilon_2^{(1)}, \varepsilon_2 \rangle = 2, 
  \\
   \nu_{\epsilon_1^{(2)}}^{(2)} (x_1) = \langle \epsilon_1^{(2)}, \varepsilon_1 \rangle = 4, \quad 
  \nu_{\epsilon_1^{(2)}}^{(2)} (x_2)  = \langle \epsilon_1^{(2)}, \varepsilon_2 \rangle = 0.
   \end{array}
\]
By Lemma \ref{lem: t2}, we obtain: 
\[
\begin{array}{c}
\nu_{\epsilon_1^{(1)}}^{(1)} (x_3) = \langle \epsilon_1^{(1)}, \g_1 \rangle = 3, \quad  
\nu_{\epsilon_2^{(1)}}^{(1)} (x_3) = \langle \epsilon_2^{(1)}, \g_1 \rangle = 1, \quad 
\nu_{\epsilon_1^{(2)}}^{(2)} (x_3) = \langle \epsilon_1^{(2)}, \g_1 \rangle = 6, \quad
\\
\nu_{\epsilon_1^{(2)}}^{(2)} (x_4) = \langle \epsilon_1^{(2)}, \g_2 \rangle = 13, \quad
\nu_{\epsilon_2^{(2)}}^{(2)} (x_4) = \langle \epsilon_2^{(2)}, \g_2 \rangle = 2.
\end{array}
\]
Finally, we apply Lemma \ref{lem: t1} to deduce that:
\[
\nu_{\epsilon_1^{(1)}}^{(1)} (x_5) = \langle \epsilon_1^{(1)}, n_1 n_2 \g_1 \rangle = 24, \quad  
\nu_{\epsilon_2^{(1)}}^{(1)} (x_5) = \langle \epsilon_2^{(1)}, n_1 n_2 \g_1 \rangle = 8, \quad 
\nu_{\epsilon_1^{(2)}}^{(2)} (x_5) = \langle \epsilon_1^{(2)}, n_2 \g_2 \rangle = 52. \quad
\]

 We have determined the following rays of the local tropicalization of 
of $\C^3$ defined by the tuple $(x_1, \dots, x_5)$: 
\begin{equation} \label{eq:tropvect}
\trop (\nu_{\epsilon_1^{(1)}}^{(1)}) = (2,0, 3,6, 24) , \quad
\trop (\nu_{\epsilon_1^{(1)}}^{(1)}) = (0, 2, 1, 2, 8), \quad
\trop (\nu_{\epsilon_1^{(1)}}^{(1)}) = (4,0,6, 13,52).
\end{equation}
The other rays of $\Theta$ are defined by the vanishing order valuations along $\X_j$, for $j =1, \dots, 5$ (see remark \ref{def:exc-div}). 
\[
\begin{array}{c}
\trop (\nu_{\epsilon_1^{(0)}}^{(0)}) = (1,0,0,0,0) , \quad
\trop( \nu_{\epsilon_2^{(0)}}^{(0)}) = (0,1,0,0,0), \quad
\trop (\nu_{\epsilon_3^{(0)}}^{(0)}) = (0,0,1, 0,0), \quad 
\\
\trop (\nu_{\epsilon_4^{(0)}}^{(0)}) = (0,0,0, 1,0), \quad 
\trop (\nu_{\epsilon_5^{(0)}}^{(0)}) = (0,0,0, 0,1).
\end{array}
\]
The local tropicalization of $Y$ in $\C^5$ 
is represented in Figure \ref{fig-LocTrop}. This figure appears as Figure 2 in the paper \cite{MMFGG}.

Remark that we can also compute the log-discrepancies of $D_i^{(\ell)}$ by using formula \eqref{eq:ld-ui}, applied to  the primitive integral vectors \eqref{eq:tropvect}.

\medskip

\setlength{\unitlength}{0.6 mm}
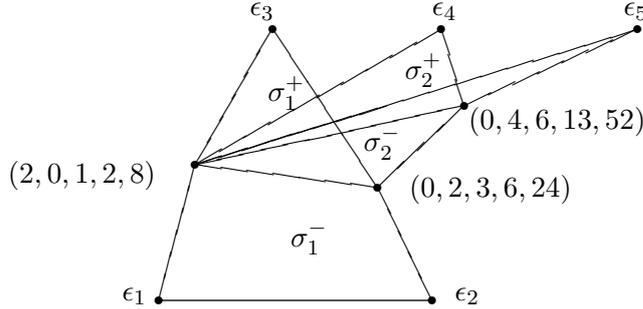
\begin{figure}[h]
\begin{center}
\begin{picture}(-60,40)(-6,-1)
\linethickness{0.1mm}
\drawline(-85,-20)(-25,-20)

\drawline(-25,-20)(-37,5)
\drawline(-85,-20)(-77,10)
\drawline(-37,5)(-60,40)
\drawline(-77,10)(-60,40)
\drawline(-77,10)(20,40)

\drawline(-77,10)(-37,5)
\drawline(-77,10)(-35,23)
\drawline(-77,10)(-23,40)
\drawline(-77,10)(-18,23)

\drawline(-18,23)(-37,5)
\drawline(-18,23)(-23,40)
\drawline(-18,23)(20,40)

\jput(-85,-20){\circle*{2}}
\jput(-25,-20){\circle*{2}}
\jput(-60,40){\circle*{2}}
\jput(-77,10){\circle*{2}}
\jput(-37,5){\circle*{2}}
\jput(-18,23){\circle*{2}}
\jput(20,40){\circle*{2}}
\jput(-23,40){\circle*{2}}

\jput(-93,-20){$\epsilon_1$}
\jput(-20,-20){$\epsilon_2$}
\jput(-65,43){$\epsilon_3$}

\jput(-118,6){$(2,0,1,2,8)$}
\jput(-30,3){$(0,2,3,6,24)$}
\jput(-56,-8){${\sigma}^-_1$}
\jput(-61,24){${\sigma}^+_1$}
\jput(-40,13){${\sigma}^-_2$}

\jput(-31,29){${\sigma}^+_2$}
\jput(17,43){$\epsilon_5$}
\jput(-25,43){$\epsilon_4$}
\jput(-17,18){$(0,4,6,13,52)$}
 \end{picture}
\end{center}

\

\

\caption[]{A sketch representing the local tropicalization of $Y$. }\label{fig-LocTrop}
\end{figure}

We apply then Theorem \ref{rem: mon-ideal}. 
We compute the values of $\xi_{\mathcal{M}}$ for those monomials 
$\mathcal{M} = x_1^{i_1} x_2^{i_2} x_3^{i_3} x_4^{i_4}$ such that $\xi_{\mathcal{M}} <1$ where 
$0 \leq i_1$, $0 \leq i_2$,  $0 \leq i_3 < 2$ and $0 \leq i_4 <4$.
We give below the
 set of jumping numbers of the multiplier ideals $\mathcal{J} (\xi D)$, for $0 < \xi <1$, 
 together with monomials which determine them. In particular, we recover the \textit{log-canonical threshold},  which is equal to $5/24$ in this example. A description of the log-canonical threshold for irreducible quasi-ordinary hypersurface singularities was given in terms of the characteristic exponents  in \cite{BGG12}.

\[
\displaystyle{
\begin{array} {|c|c|c|c|c|c|c|c|c|c|c|}
\hline
\mathcal{M} & 1 & x_1 & y& x_1^2 & x_1^3 & x_1y & x_1^3x_2  & z_1 & x_1^2y & x_1^3y
\\ \hline 
 \xi_{\mathcal{M}} & \frac{5}{24} & \frac{15}{52} & \frac{17}{52} & \frac{19}{52} & \frac{3}{8} 
 & \frac{21}{52} & \frac{23}{52} & \frac{11}{24} & \frac{25}{52}& \frac{4}{8}
 \\  \hline \hline 
 \mathcal{M} & x_1^4x_2 &x_1z_1 & x_1^3x_2y& yz_1 & x_1^5x_2 &x_1^2z_1 &x_1^2z_1 & x_1^4x_2y & x_1yz_1 & x_1^6x_2^2 
 \\ \hline
\xi_\mathcal{M}  & \frac{27}{52} & \frac{28}{52}& \frac{29}{52}&  \frac{30}{52} & \frac{31}{52} & \frac{32}{52}& \frac{5}{8} &\frac{33}{52} &\frac{34}{52} &  \frac{35}{52} 
  \\ \hline \hline 
 \mathcal{M} &x_1^3x_2z_1 & z_1^2 & x_1^5x_2y & x_1^2yz_1 &x_1^3yz_1 & x_1^4x_2z_1 & x_1z_1^2& x_1^3x_2yz_1 & yz_1^2 & x_1^5x_2z_1
   \\ \hline 
\xi_\mathcal{M} & \frac{36}{52} & \frac{17}{24} & \frac{37}{52}& \frac{38}{52} & \frac{39}{52} &\frac{40}{52} & \frac{41}{52} & \frac{42}{52} & \frac{43}{52} & \frac{44}{52}
   \\ \hline \hline 
\mathcal{M} & x_1^2z_1^2 & x_1^3z_1^2 & x_1^4x_2yz_1 & x_1yz_1^2 & x_1^6x_2^2z_1 & x_1^3x_2z_1^2 & z_1^3 & x_1^5x_2yz_1 & x_1^2yz_1^2 & 
   \\ \hline
\xi_\mathcal{M} & \frac{45}{52} & \frac{7}{8} & \frac{46}{52} & \frac{47}{52} & \frac{48}{52} & \frac{49}{52} &\frac{23}{24}& \frac{50}{52} &\frac{51}{52} & 
   \\ \hline  
   \end{array}
   }
\]
\end{example}

%
%%  \nocite*
%  \bibliographystyle{amsalpha}
% 
% \renewcommand{\MR}[1]{}
%
%\bibliography{referencias-qo}
%
%
%\end{document}

% \bibliographystyle{abbrv}
% \bibliography{referencias-qo.bib}

 \end{document}